\documentclass[12pt,twoside,sort&compress ]{article} 
\usepackage{mathrsfs,amsfonts,amsmath,amsbsy,tikz}\usepackage{fontenc}\usepackage{textcomp}
\usepackage{graphicx}
\usepackage{amssymb}
\usepackage{graphicx}

\def\le{\leq}

\def\bbb{\begin{eqnarray*}}

\def\eee{\end{eqnarray*}}

\begin{document} 
\baselineskip=18pt
\begin{center}
\vspace{-0.6in}{\large \bf Inequalities among eigenvalues of different \\[0.1in]
self-adjoint discrete  Sturm-Liouville  problems}
\\ [0.2in]
 Hao Zhu,\ \ Yuming Shi $^{\dag}$  \\
\vspace{0.15in}  Department of Mathematics, Shandong University\\
Jinan, Shandong 250100, P. R. China\\

\footnote{$^\dag$ The corresponding author}
\footnote{ \;\;Email addresses: haozhusdu@163.com(H. Zhu), ymshi@sdu.edu.cn(Y. Shi).}
\footnote{ \;\;This research is supported by the NNSF of China (Grants 11571202).}
\end{center}

{\bf Abstract.}  In this paper, inequalities
 among eigenvalues of different self-adjoint discrete Sturm-Liouville problems  are established.
 For a fixed discrete Sturm-Liouville equation, inequalities among eigenvalues for different  boundary conditions  are given. For a fixed  boundary condition, inequalities among eigenvalues for different equations are given. These results are obtained by applying continuity and discontinuity of the $n$-th eigenvalue function,
  monotonicity in some direction of the $n$-th eigenvalue function, which were given in our previous papers, and natural loops in the space of boundary conditions. Some results  generalize the relevant existing results about inequalities among eigenvalues of different Sturm-Liouville problems.

\medskip

\medskip

{2010 {\bf \it Mathematics Subject Classification}}: 39A70; 34B24; 15A42; 47A75.

\bigskip

\noindent{\bf 1. Introduction}\medskip

A self-adjoint discrete Sturm-Liouville problem (briefly, SLP) considered in the present paper consists of a symmetric discrete Sturm-Liouville equation (briefly, SLE) \vspace{-0.2cm}
$$ -\nabla(f_{n}\Delta y_{n})+q_{n}y_{n}=\lambda w_{n}y_{n}, \;\;\;\; \ \  n\in[1,N],                                                                 \eqno(1.1)
\vspace{-0.2cm} $$
 and a self-adjoint boundary condition (briefly, BC)\vspace{-0.2cm}
$$
A \left (\begin{array} {cc}y_{0}\\
f_{0}\triangle y_{0}\end{array}   \right )
+B\left ( \begin{array} {cc}y_{N}\\
f_{N}\triangle y_{N}
\end{array}   \right )=0,                                                                                                                             \eqno(1.2)
\vspace{-0.2cm}
$$
where $N \geq 2$ is an integer, $\Delta $ and $\nabla$ are the forward and backward difference
operators, respectively, i.e., $\Delta y_n=y_{n+1}-y_n$ and $\nabla y_n=y_n-y_{n-1}$;
$f=\{f_n\}_{n=0}^{N}$, $q=\{q_n\}_{n=1}^{N}$ and $w=\{w_n\}_{n=1}^{N}$ are real-valued sequences such that
\vspace{-0.2cm}$$f_n\neq0 \;{\rm for}\;\; n\in [0,N],
\;\;w_n>0\;{\rm for}\;\; n\in [1,N];                                                                                                                 \eqno(1.3)
\vspace{-0.2cm}$$
$\lambda$ is the spectral parameter; the interval $[M,N]$ denotes the set of integers
 $\{M,M+1,\cdots,N\}$; $A$ and $B$ are
$2\times2$ complex matrices such that rank(A, B)=2, and satisfy the following  self-adjoint boundary  condition:
\vspace{-0.05cm} $$
 A\left (\begin{array} {cc} 0&1\\
-1&0\end{array}   \right )A^{*}=
B\left (\begin{array} {cc} 0&1\\
-1&0\end{array}   \right )B^{*},\vspace{-0.05cm}                                                                                                       \eqno(1.4)
$$
where $A^{*}$ denotes the complex conjugate transpose of $A$.

Throughout this paper, by $\mathbb{C}$, $\mathbb{R}$, and $\mathbb{Z}$ denote the sets of the complex
numbers, real numbers, and integer numbers, respectively; and by $\bar{z}$ denote the
complex conjugate of $z\in \mathbb{C}$. Moreover, when a capital Latin letter stands for a
matrix, the entries of the matrix are denoted by the corresponding
lower case letter with two indices. For example, the entries of a
matrix $C$ are $c_{ij}$'s.

As it is mentioned in [1, 7], the discrete SLP (1.1)--(1.2) can be applied to many fields, ranging from mechanics, to network theory, and to probability theory. The eigenvalues of (1.1)--(1.2) play an important role in studying these physical problems and they change as the SLP changes. Thus, it is naturally important to compare the
eigenvalues of different SLPs. In this paper, we shall establish inequalities among eigenvalues of different SLPs.

Recall that a self-adjoint continuous SLP consists of a differential SLE
$$ -(p(t)y')'+q(t)y=\lambda w(t)y, \;\;\ \  t\in(a,b),
                                                                                                                                                        \eqno(1.5)
 \vspace{-0.2cm} $$
 and a  BC \vspace{-0.2cm}
$$
A \left (\begin{array} {cc}y(a)\\
(py')(a)\end{array}   \right )
+B\left ( \begin{array} {cc}y(b)\\
(py')(b)
\end{array}   \right )=0,
                                                                                                                                                       \eqno(1.6)
\vspace{-0.05cm}$$
where $-\infty< a<b<+\infty$; $1/p,q,w\in L((a,b),\mathbb{R}),$ $p,w>0$ almost everywhere  on $(a,b)$, while $L((a,b),\mathbb{R})$ denotes the space of Lebesgue integrable
real functions on $(a,b)$; $A$ and $B$ are
$2\times2$ complex matrices such that rank(A, B)=2 and (1.4) holds.
For a fixed differential SLE (1.5), inequalities among eigenvalues for different  self-adjoint  BCs  have been extensively studied by many authors (cf., e.g. [3--6, 8, 10, 12--15, 20, 21]).  Using the variational method,  Courant and Hilbert in [5] gave inequalities among eigenvalues for different separated BCs.
Using the Pr\"{u}fer transformation of (1.5),  Coddington and Levinson in [4] gave the classical  inequalities among eigenvalues for periodic ,  antiperiodic , Dirichlet and Neumann BCs   under some conditions on the coefficients of (1.5).
See also [20]. For  an arbitrary coupled self-adjoint BC, Eastham and his coauthors in [6, Theorem 3.2] identified two separated
BCs corresponding to the Dirichlet and Neumann BCs in the
above case, and established analogous inequalities.
Their proof also depends on the Pr\"{u}fer transformation of (1.5). See also [8].  These inequalities are extended to singular SLPs and other cases [3, 10, 12--14]. Using natural loops in the space of self-adjoint  BCs, Peng and his coauthors in [15]  gave a short proof of [6, Theorem 3.2],  and obtained  new general inequalities. See Theorem 4.53 in [15].

Next,  we recall the related existing results of inequalities among eigenvalues of different self-adjoint discrete SLPs (1.1)--(1.2). For a fixed self-adjoint BC (1.2),
inequalities among eigenvalues for different coefficients of (1.1) were obtained by   Rayleigh¡¯s
principles and minimax theorems  in [16]. Then these results were extended to higher-order discrete vector SLPs in [17].
For a fixed equation (1.1), by using
some oscillation results obtained in [1] and some spectral results of (1.1)--(1.2) obtained in [16],
inequalities among eigenvalues for periodic,  antiperiodic, and Dirichlet BCs were given  under the assumption that $f_n>0$, $1\leq n\leq N-1$, and $f_0=f_N=1$ in [19]. Under the same conditions of the coefficients of (1.1) and by a similar method used in [19], the above results in [19] were extended to a class of coupled BCs  in [18].

The aim of the present paper is to establish more general inequalities among eigenvalues of different SLPs (1.1)--(1.2).
For a fixed equation, inequalities among eigenvalues for different separated BCs are established in Theorem 3.1, and  among eigenvalues for different BCs
in a natural loop are established in Theorems 3.2--3.5. Then, the inequalities in Theorems 3.2--3.5 are applied to compare eigenvalues for  coupled
BCs with those for some certain separated ones (see Theorems 3.6--3.11), and  eigenvalues for different coupled BCs (see Theorem 3.12). The inequalities in Theorems 3.6--3.11 extend those in [18, Theorem 3.1] to a more general case. For a fixed BC, inequalities among eigenvalues for equations with different coefficients and weight functions are established in Theorem 4.1, which generalize those in [16, Theorem 5.5] and [17, Theorem 3.6] in the second-order case.
Combining the above results, one can establish inequalities among eigenvalues of SLPs with different equations and BCs (see Corollary 4.1 and Remark 4.2).

The method used in the present paper is different from those used in [4, 6, 8, 18--20]. On the one hand,  the approaches used in [4, 6, 8, 20]  in the continuous case depend on the Pr\"{u}fer transformation of (1.5). Although  the  Pr\"{u}fer transformation in discrete version were given in [2],  some of its properties in continuous version can not be extended to the discrete one and thus
similar methods used in [4, 6, 8, 20] are difficultly employed in studying the discrete problem.
On the other hand, the variational method used in [16, 17] is restricted to compare eigenvalues, which have the same index, of different SLPs, and  it seems to us that the method used in [18, 19] in the discrete case is hardly extended to a more general case.
Note that  the method used in [15] does not depend on the Pr\"{u}fer transformation of (1.5). So the similarity of the self-adjoint
 BCs for continuous and discrete SLPs [11, 23] makes it possible and convenient to generalize a similar approach used in [15] to the discrete case. There are three major ingredients, which will be used by  this method: (1) the continuity and discontinuity of the $n$-th eigenvalue function, which were studied in [22];
 (2) the monotonicity  of the $n$-th eigenvalue function, which can be deduced from [22, 23]; (3) natural loops in the space of self-adjoint BCs, which can be obtained in a similar way as  that in [15].
Thus, the work in the present paper, to some extent, can be regarded as a discrete analog of that in [15] and a continuation of our present works [22, 23].

 This paper is organized as follows.  Section 2  gives some preliminaries. Some notations are introduced and some lemmas are recalled. Especially, natural loops in space of self-adjoint BCs,  are presented. In Section 3, inequalities among eigenvalues for different boundary conditions are given.
In Section 4, inequalities among eigenvalues for different equations are established.
\bigskip

\noindent{\bf 2. Preliminaries }\medskip

In this section, some notations and lemmas are introduced. This section is divided into two parts.  In Section 2.1, topology on space of SLPs
 and several useful properties of eigenvalues are recalled. In Section 2.2, natural loops in space of self-adjoint BCs established in [15] are presented. \medskip

\noindent{\bf 2.1. Space of SLPs and properties of eigenvalues }\medskip

Let the SLE (1.1) be abbreviated as $(1/f,q,w)$.  Then the space of the SLEs can be written as \vspace{-0.1cm}$$\begin{array}{rrll}\Omega_N^{\mathbb{R},+} :=\{(1/f,q,w)\in\mathbb{R}^{3N+1}: {\rm (1.3)\; holds}\},\end{array} \vspace{-0.2cm}$$
and is equipped with the topology deduced from the real space $ \mathbb{R}^{3N+1}$. Note that $\Omega_N^{\mathbb{R},+}$ has $2^{N+1}$ connected components.
Bold faced lower case Greek letters, such as $ \pmb\omega$, are used to denote elements of $\Omega_N^{\mathbb{R},+}$.

The quotient space
\vspace{-0.2cm}
$$\mathcal A^{\mathbb{C}} :=\raise2pt\hbox{${\rm M}^*_{2,4}(\mathbb C)$}/\lower3pt\hbox{${\rm
GL}(2,\mathbb C)$},
\vspace{-0.2cm}$$
equipped with the quotient topology, is taken as the space of general BCs; that is, each BC is an equivalence class of coefficient matrices of
system (1.2), where
\vspace{-0.1cm}$$\begin{array}{l}M_{2,4}^{*}(\mathbb{C}):=\{2\times4 \;{\rm complex\; matrix\; (A,B):}\;{\rm rank}(A,B)=2\},\\[1.0ex]
GL(2,\mathbb{C}):=\{2\times2\; {\rm comlplex \;matrix}\; T:{\rm det}\; T \neq0\}.\end{array}\vspace{-0.1cm} $$
The BC represented by system (1.2) is denoted by $[A\,|\,B]$. Bold faced capital Latin letters, such as $\mathbf{A}$, are also used for BCs.
The space of self-adjoint BCs is denoted by $\mathcal{B}^\mathbb{C}$. The following result gives the topology and geometric structure of $\mathcal{B}^\mathbb{C}$.\medskip

\noindent{\bf Lemma 2.1} {\rm [23, Theorem 2.2]}. {\it The space $\mathcal{B}^\mathbb{C}$ equals the union of the following relative open sets:\vspace{-0.2cm}
$$\begin{array} {cccc} \mathcal{O}_{1,3}^{\mathbb{C}}=\left \{\left [\begin{array} {cccc}1&a_{12}&0&\bar{z}\\
0&z&-1&b_{22}\end{array}  \right ]:\; a_{12},b_{22}\in\mathbb{R},z\in \mathbb{C}\right \},\vspace{3mm}\\
 \mathcal{O}_{1,4}^{\mathbb{C}}=\left \{\left [\begin{array} {cccc}1&a_{12}&\bar{z}&0\\
0&z&b_{21}&1\end{array}  \right ]:\; a_{12},b_{21}\in\mathbb{R},z\in \mathbb{C}\right \},\vspace{3mm}\\
\mathcal{O}_{2,3}^{\mathbb{C}}=\left \{\left [\begin{array} {cccc}a_{11}&-1&0&\bar{z}\\
z&0&-1&b_{22}\end{array}  \right ]:\; a_{11},b_{22}\in\mathbb{R},z\in \mathbb{C}\right \},\vspace{3mm}\\
 \mathcal{O}_{2,4}^{\mathbb{C}}=\left \{\left [\begin{array} {cccc}a_{11}&-1&\bar{z}&0\\
z&0&b_{21}&1\end{array}  \right ]:\; a_{11},b_{21}\in\mathbb{R},z\in \mathbb{C}\right \}.\end{array}\vspace{-0.05cm}
                                                                                                                                                      \eqno(2.1)
$$
Moreover, $\mathcal{B}^\mathbb{C}$ is a connected and compact real-analytic manifold of dimension {\rm 4}.}\medskip

 Lemma 2.1 says that $\mathcal{O}_{1,3}^{\mathbb{C}}$, $\mathcal{O}_{1,4}^{\mathbb{C}}$, $\mathcal{O}_{2,3}^{\mathbb{C}}$, and $\mathcal{O}_{2,4}^{\mathbb{C}}$ together form an atlas of local coordinate systems on $\mathcal{B}^\mathbb{C}$.

          The space $\Omega_N^{\mathbb{R},+}\times\mathcal{B}^\mathbb{C}$ of the SLPs is a real-analytic manifold of dimension
$3N+5$ and has $2^{N+1}$ connected components.

      The following result gives the canonical forms of separated and coupled self-adjoint BCs.\medskip

\noindent{{\bf Lemma 2.2} {\rm [21, Theorem 10.4.3]}. {\it A separated self-adjoint BC can be written as \vspace{-0.1cm}
$$\mathbf{S}_{\alpha,\beta}:=\left [\begin{array} {cccc}\cos\alpha&-\sin\alpha&0&0\\
0&0&\cos\beta&-\sin\beta\end{array}  \right ],                                                                                                        \eqno(2.2)
\vspace{-0.2cm}$$
where $\alpha\in[0,\pi),\beta\in (0,\pi];$
and a coupled self-adjoint BC can be written as \vspace{-0.2cm}
$$[e^{i\gamma}K\,|\,-I],\vspace{-0.2cm}$$
where $\gamma\in(-\pi,\pi],\;\;K\in SL(2,\mathbb{R}):=\{2\times2 \;real\; matrix\; M:{\rm det}M=1\}.$}\medskip

In particular, $\mathbf{S}_{0,\pi}$ is called the Dirichlet BC; $\mathbf{S}_{0,\beta}$ for any $\beta\in(0,\pi]$ or $\mathbf{S}_{\alpha,\pi}$ for any $\alpha\in[0,\pi)$ is called the BC which is Dirichlet at an endpoint.   By $\mathcal{B}_S$ and $\mathcal{B}_C$ denote the space of separated self-adjoint BCs and that of coupled self-adjoint BCs, respectively. Then $\mathcal{B}^{\mathbb{C}}=\mathcal{B}_S\cup\mathcal{B}_C$, and $\mathcal{B}_C$ is an open set of
 $\mathcal{B}^{\mathbb{C}}$.

Next, several properties of eigenvalues are presented.
For each $\lambda\in\mathbb{C}$, let $\phi(\lambda)=\{\phi_n(\lambda)\}_{n=0}^{N}$ and $\psi(\lambda)=\{\psi_n(\lambda)\}_{n=0}^{N}$ be the solutions of (1.1)
satisfying the following initial conditions: \vspace{-0.2cm}
$$
\phi_0(\lambda)=1, f_0\Delta\phi_0(\lambda)=0; \;\;
\psi_0(\lambda)=0, \ f_0\Delta\psi_0(\lambda)=1.
\vspace{-0.2cm}$$
Then the leading terms of $\phi_N (\lambda)$, $\psi_N (\lambda)$, $f_N\Delta\phi_N (\lambda)$,
and $f_N\Delta\psi_N (\lambda)$ as polynomials of $\lambda$ are \vspace{-0.2cm} $$\begin{array}{cccc}(-1)^{N-1}\left(\prod\limits_{i=1}^{N-1}({w_i}/{f_i})\right)
\lambda^{N-1},&(-1)^{N-1}\left(({1}/{f_0})\prod\limits_{i=1}^{N-1}({w_i}/{f_i})\right)
\lambda^{N-1},
\\[2.0ex](-1)^{N}\left(w_N\prod\limits_{i=1}^{N-1}({w_i}/{f_i})\right)
\lambda^{N},&(-1)^{N}\left(({w_N}/{f_0})\prod\limits_{i=1}^{N-1}({w_i}/{f_i})\right)
\lambda^{N},\end{array}                                                                                                                               \eqno(2.3)
\vspace{-0.2cm} $$
respectively. See [23] for details.
\medskip

The following result says that  the eigenvalues
of a given SLP can be determined by a polynomial.\medskip

\noindent{\bf Lemma 2.3 } {\rm [23, Lemma 3.2 and Lemma 3.3]}. {\it A number $\lambda \in
\mathbb{C}$ is an eigenvalue of each given SLP
{\rm(1.1)}--{\rm(1.2)}
 if and only if $\lambda$ is a zero of the polynomial\vspace{-0.2cm}
$$\Gamma(\lambda) =\det A +\det B +G(\lambda),\vspace{-0.2cm}$$
where}\vspace{-0.2cm}
$$G(\lambda):=c_{11}\phi_N(\lambda)+c_{12}\psi_N(\lambda)+c_{21}f_N\Delta\phi_N(\lambda)
+c_{22}f_N\Delta\psi_N(\lambda),\vspace{-0.2cm}$$
\vspace{-0.2cm}$$\begin{array}{rrll}c_{11}:=a_{22}b_{11}-a_{12}b_{21},&c_{12}:=a_{11}b_{21}-a_{21}b_{11},\\[1.0ex]
c_{21}:=a_{22}b_{12}-a_{12}b_{22},&c_{22}:=a_{11}b_{22}-a_{21}b_{12}.\end{array} $$

 Let $(\pmb\omega,\mathbf{A})\in\Omega_N^{\mathbb{R},+}\times\mathcal{B}^{\mathbb{C}}$. Set \vspace{-0.2cm}
$$r=r(\pmb\omega,\mathbf{A}):={\rm rank} \left(\begin{array}{cc}-a_{11}+f_0a_{12}&b_{12}\\-a_{21}+f_0a_{22}&b_{22}\end{array}\right).\eqno(2.4)\vspace{-0.2cm}
$$
 Obviously, $0\leq r \leq 2$.
The following result
establishes the relationship between analytic and geometric  multiplicities
of each eigenvalue of a given SLP and gives a formula for counting the number of eigenvalues.
\medskip

\noindent{\bf Lemma 2.4 } [23, Lemma 3.4 and Theorem 3.3]. {\it For any fixed self-adjoint SLP {\rm(1.1)--(1.2)}, all its eigenvalues are real,
the number of its eigenvalues   is equal to  $N-2+r$,  where $r$ is defined by {\rm(2.4)}, and the analytic and geometric multiplicities of each of its eigenvalue   are the same. }\medskip

Lemma 2.4  can also be deduced from [16, Theorem 4.1] and [17, Theorem 4.3]. By Lemma 2.4, we shall only say the multiplicity of an eigenvalue without specifying its analytic and geometric multiplicities.  Based on these results, the problem (1.1)--(1.2) has $k=N-2+r$ eigenvalues (counting multiplicities), which  can be arranged in the following non-decreasing order:
\vspace{-0.2cm}$$\lambda_0\le \lambda_1\le \lambda_2 \le \cdots\le\lambda_{k-1}.\vspace{-0.2cm}$$
The $n$-th eigenvalue $\lambda_n$ can be considered as a function in the space of the SLPs, called the $n$-th eigenvalue function. The following result gives a necessary and sufficient condition for all the eigenvalue functions to be continuous in a  set of space of SLPs.\medskip

 \noindent{\bf Lemma 2.5}  [22, Theorem 2.1]. {\it Let $\mathcal{O}$ be a  set of $\Omega_N^{\mathbb{R},+}\times\mathcal{B}^{\mathbb{C}}$.
Then the number of eigenvalues of each $(\pmb\omega,\mathbf {A})\in\mathcal{O}$ is equal if and only if all  the eigenvalue functions restricted in $\mathcal{O}$  are
 continuous. Furthermore, if  $\mathcal{O}$ is a connected set of $\Omega_N^{\mathbb{R},+}\times\mathcal{B}^{\mathbb{C}}$, then each eigenvalue function is locally
a continuous eigenvalue branch in $\mathcal{O}$.}

\bigskip

\noindent{\bf 2.2. Natural loops in the space of self-adjoint boundary conditions  }\medskip

In this subsection, natural loops in the space of self-adjoint BCs are presented. We remark that these natural loops will  play an important role in studying inequalities among eigenvalues for coupled BCs and those for  some certain separated ones.
\medskip

\noindent{\bf Lemma 2.6} [15, Lemma 3.1].  {\it In $\mathcal{B}^{\mathbb{C}}$, we have the following limits:
$$\vspace{-0.2cm}\begin{array} {cccc}\mathbf{S}_1:=\lim\limits_{s\rightarrow\pm\infty}\left [\begin{array} {cccc}1&s&\bar{z}&0\\
0&z&b_{21}&1\end{array}  \right ]=\left [\begin{array} {cccc}0&1&0&0\\
0&0&b_{21}&1\end{array}  \right ],\vspace{3mm}\\
\mathbf{S}_2:=\lim\limits_{t\rightarrow\pm\infty}\left [\begin{array} {cccc}1&a_{12}&\bar{z}&0\\
0&z&t&1\end{array}  \right ]=\left [\begin{array} {cccc}1&a_{12}&0&0\\
0&0&1&0\end{array}  \right ],\vspace{3mm}\\
\mathbf{S}_3:=\lim\limits_{s\rightarrow\pm\infty}\left [\begin{array} {cccc}s&-1&\bar{z}&0\\
z&0&b_{21}&1\end{array}  \right ]=\left [\begin{array} {cccc}1&0&0&0\\
0&0&b_{21}&1\end{array}  \right ],\vspace{3mm}\\
\mathbf{S}_4:=\lim\limits_{t\rightarrow\pm\infty}\left [\begin{array} {cccc}a_{11}&-1&\bar{z}&0\\
z&0&t&1\end{array}  \right ]=\left [\begin{array} {cccc}a_{11}&-1&0&0\\
0&0&1&0\end{array}  \right ],\vspace{3mm}\\
\mathbf{S}_5:=\lim\limits_{s\rightarrow\pm\infty}\left [\begin{array} {cccc}s&-1&0&\bar{z}\\
z&0&-1&b_{22}\end{array}  \right ]=\left [\begin{array} {cccc}1&0&0&0\\
0&0&-1&b_{22}\end{array}  \right ],\vspace{3mm}\\
\mathbf{S}_6:=\lim\limits_{t\rightarrow\pm\infty}\left [\begin{array} {cccc}a_{11}&-1&0&\bar{z}\\
z&0&-1&t\end{array}  \right ]=\left [\begin{array} {cccc}a_{11}&-1&0&0\\
0&0&0&1\end{array}  \right ],\vspace{3mm}\\
\mathbf{S}_7:=\lim\limits_{s\rightarrow\pm\infty}\left [\begin{array} {cccc}1&s&0&\bar{z}\\
0&z&-1&b_{22}\end{array}  \right ]=\left [\begin{array} {cccc}0&1&0&0\\
0&0&-1&b_{22}\end{array}  \right ],\vspace{3mm}\\
\mathbf{S}_8:=\lim\limits_{t\rightarrow\pm\infty}\left [\begin{array} {cccc}1&a_{12}&0&\bar{z}\\
0&z&-1&t\end{array}  \right ]=\left [\begin{array} {cccc}1&a_{12}&0&0\\
0&0&0&1\end{array}  \right ].
\end{array}$$}

The following result  gives some natural loops  in  $\mathcal{B}^{\mathbb{C}}$.\medskip

\noindent{\bf Lemma 2.7} [15, Lemma 3.7]. {\it  \begin{itemize}\vspace{-0.2cm}
\item[{\rm (i)}]  Every   BC $\mathbf{A}\in\mathcal{O}_{1,4}^{\mathbb{C}}$ lies on the following two simple real-analytic loops in
$\mathcal{B}^{\mathbb{C}}$:
$$\vspace{-0.2cm}\begin{array} {cccc}\mathcal{C}_{1,4,z, b_{21}}=\left\{ \mathbf{A}(s):=\left [\begin{array} {cccc}1&s&\bar{z}&0\\
0&z&b_{21}&1\end{array}  \right ],\;s\in\mathbb{R}\right\}\cup\left\{\mathbf{S}_1\right\},
\\[2.0ex]
\mathcal{C}_{1,4,z, a_{12}}=\left\{ \hat{\mathbf{A}}(t):=\left [\begin{array} {cccc}1&a_{12}&\bar{z}&0\\
0&z&t&1\end{array}  \right ],\;t\in\mathbb{R}\right\}\cup\left\{\mathbf{S}_2\right\}.\end{array}\vspace{-0.15cm}\vspace{-0.2cm}$$
\vspace{-0.2cm}
\item[{\rm (ii)}]  Every   BC $\mathbf{A}\in\mathcal{O}_{2,4}^{\mathbb{C}}$ lies on the following two simple real-analytic loops in
$\mathcal{B}^{\mathbb{C}}$:
$$\vspace{-0.2cm}\begin{array} {cccc}\mathcal{C}_{2,4,z, b_{21}}=\left\{ \mathbf{B}(s):=\left [\begin{array} {cccc}s&-1&\bar{z}&0\\
z&0&b_{21}&1\end{array}  \right ],\;s\in\mathbb{R}\right\}\cup\left\{\mathbf{S}_3\right\},
\\[2.0ex]
\mathcal{C}_{2,4,z, a_{11}}=\left\{ \hat{\mathbf{B}}(t):=\left [\begin{array} {cccc}a_{11}&-1&\bar{z}&0\\
z&0&t&1\end{array}  \right ],\;t\in\mathbb{R}\right\}\cup\left\{\mathbf{S}_4\right\}.\end{array}\vspace{-0.1cm}$$
\item[{\rm (iii)}]  Every  BC $\mathbf{A}\in\mathcal{O}_{2,3}^{\mathbb{C}}$ lies on the following two simple real-analytic loops in
$\mathcal{B}^{\mathbb{C}}$:
$$\vspace{-0.2cm}\begin{array} {cccc}\mathcal{C}_{2,3,z, b_{22}}=\left\{\mathbf{C}(s):=\left [\begin{array} {cccc}s&-1&0&\bar{z}\\
z&0&-1&b_{22}\end{array}  \right ],\;s\in\mathbb{R}\right\}\cup\left\{\mathbf{S}_5\right\},
\\[2.0ex]
$$\vspace{-0.2cm}\mathcal{C}_{2,3,z, a_{11}}=\left\{ \hat{\mathbf{C}}(t):=\left [\begin{array} {cccc}a_{11}&-1&0&\bar{z}\\
z&0&-1&t\end{array}  \right ],\;t\in\mathbb{R}\right\}\cup\left\{\mathbf{S}_6\right\}.\end{array}\vspace{-0.2cm}$$
\vspace{-0.2cm}
\item[{\rm (iv)}] Every   BC $\mathbf{A}\in\mathcal{O}_{1,3}^{\mathbb{C}}$ lies on the following two simple real-analytic loops in
$\mathcal{B}^{\mathbb{C}}$:
$$\vspace{-0.2cm}\begin{array} {cccc}\mathcal{C}_{1,3,z, b_{22}}=\left\{\mathbf{D}(t):= \left [\begin{array} {cccc}1&s&0&\bar{z}\\
0&z&-1&b_{22}\end{array}  \right ],\;s\in\mathbb{R}\right\}\cup\left\{\mathbf{S}_7\right\},
 \\[2ex]
\mathcal{C}_{1,3,z, a_{12}}=\left\{\hat{\mathbf{D}}(t):= \left [\begin{array} {cccc}1&a_{12}&\bar{z}&0\\
0&z&-1&t\end{array}  \right ],\;t\in\mathbb{R}\right\}\cup\left\{\mathbf{S}_8\right\}.\end{array}\vspace{-0.2cm}$$\vspace{-0.2cm}
 \end{itemize}}

 \noindent{\bf Remark 2.1.} $\mathcal{C}_{1,4,z, b_{21}}\backslash \{\mathbf{A}\}$ is connected for any fixed $ \mathbf{A}\in\mathcal{C}_{1,4,z, b_{21}}$. Similar result holds for other natural loops in (i)--(iv) of Lemma 2.7. For each $1\leq i\leq 8$, $\mathbf{S}_i$ is called a limit boundary condition (briefly, LBC) in the corresponding  natural loop. Note that all the LBCs are separated ones.
\bigskip

\noindent{\bf 3. Inequalities among eigenvalues for different boundary conditions }\medskip

In this section, for any fixed equation, inequalities among eigenvalues for different  BCs are established. This section is divided into four parts.   In Subsections 3.1--3.4, inequalities among eigenvalues  for different separated BCs, among eigenvalues  for different  BCs in a natural loop,
 among eigenvalues for coupled BCs and those for  some certain separated ones, and among eigenvalues for different coupled BCs are established, respectively. \medskip

\noindent{\bf 3.1. Inequalities among eigenvalues for separated BCs} \medskip

In this subsection, we shall first compare the eigenvalues for different separated BCs $\mathbf{S}_{\alpha,\beta}$ in two directions $\alpha$ and $\beta$, separately. Then we  give an application to
compare eigenvalues for an arbitrary separated BC with those for the BCs which are Dirichlet at an endpoint.

For convenience,
 denote $\lambda_n(\alpha,\beta):=\lambda_n(\mathbf{S}_{\alpha,\beta})$ for short, and
$\xi:=\arctan(-1/ f_0)+\pi$ if $ f_0>0$; $\xi:=\arctan(-1/ f_0)$ if $ f_0<0$.   \medskip

\noindent{{\bf Theorem  3.1.} {\it Fix a difference equation $\pmb\omega=(1/f,q,w)$. Then $(\pmb\omega,\mathbf{S}_{\alpha,\beta})$ has exactly $N$ eigenvalues if $\alpha\neq\xi$ and $\beta\neq\pi$;  exactly $N-1$ eigenvalues if either $\alpha\neq\xi$ and $\beta=\pi$ or $\alpha=\xi$ and $\beta\neq\pi$; and exactly $N-2$ eigenvalues if $\alpha=\xi$ and $\beta=\pi$. Further, for any  $0\leq\alpha_1<\alpha_2<\xi\leq\alpha_3<\alpha_4<\pi$, and $0<\beta_1<\beta_2\leq\pi$, we have that
\begin{itemize}\vspace{-0.2cm}
\item[{\rm (i)}]  the eigenvalues of the SLPs $(\pmb\omega,\mathbf{S}_{\alpha_i,\beta_0})$ for any $\beta_0\in(0,\pi)$, $i=1,\cdots,4$, satisfy the following inequalities:
 \vspace{-0.2cm}$$\begin{array} {cccc}\lambda_0(\alpha_2,\beta_0)<\lambda_0(\alpha_1,\beta_0)<\lambda_0(\alpha_4,\beta_0)<\lambda_0(\alpha_3,\beta_0)<\lambda_1(\alpha_2,\beta_0)<
 \lambda_1(\alpha_1,\beta_0)\\[1.0ex]
 <\lambda_1(\alpha_4,\beta_0)<\lambda_1(\alpha_3,\beta_0)<\cdots
 <
\lambda_{N-2}(\alpha_2,\beta_0) <\lambda_{N-2}(\alpha_1,\beta_0)<\\[1.0ex]\lambda_{N-2}(\alpha_4,\beta_0)<\lambda_{N-2}(\alpha_3,\beta_0)<\lambda_{N-1}(\alpha_2,\beta_0)<\lambda_{N-1}(\alpha_1,\beta_0)
<\lambda_{N-1}(\alpha_4,\beta_0),\end{array}\vspace{-0.2cm}$$
 and in addition, $\lambda_{N-1}(\alpha_4,\beta_0)<\lambda_{N-1}(\alpha_3,\beta_0)$ if $\alpha_3\neq\xi$;
\item[{\rm (ii)}]  similar results in {\rm(i)}  hold
   with $N-2$ and $N-1$ replaced by  $N-3$ and $N-2$, respectively, in the case that  $\beta_0=\pi$;
\item[{\rm (iii)}]  the eigenvalues of the SLPs $(\pmb\omega, \mathbf{S}_{\alpha_0,\beta_j})$ for any $\alpha_0\in[0,\xi)\cup(\xi,\pi)$, $j=1,2$, satisfy the following inequalities:
\vspace{-0.2cm}$$\begin{array} {llll}
\lambda_0(\alpha_0,\beta_1)<\lambda_0(\alpha_0,\beta_2)<\lambda_1(\alpha_0,\beta_1)<\lambda_1(\alpha_0,\beta_2)<\cdots
 <\\[1.0ex]
 <\lambda_{N-2}(\alpha_0,\beta_1)<\lambda_{N-2}(\alpha_0,\beta_2)<\lambda_{N-1}(\alpha_0,\beta_1),\end{array}\vspace{-0.2cm}$$
 and in addition, $\lambda_{N-1}(\alpha_0,\beta_1)<\lambda_{N-1}(\alpha_0,\beta_2)$ if $\beta_2\neq\pi$;
\item[{\rm (iv)}]  similar results in {\rm(iii)}  hold
   with $N-2$ and $N-1$ replaced by $N-3$ and $N-2$, respectively, in the case that $\alpha_0=\xi$.
\end{itemize}}

\noindent{\bf Proof.}
The number of eigenvalues of $(\pmb\omega,\mathbf{A})$ in each case can be obtained by Lemma 2.4. Firstly, we show that (i) holds. Let $\beta_0\in(0,\pi)$.
By (i) of Corollary 4.2 in [22],
the $n$-th eigenvalue functions $\lambda_n(\alpha,\beta_0)$  are strictly decreasing  in $\alpha\in[0,\xi)$ or $\alpha\in(\xi,\pi)$
 for all $0\leq n\leq N-1$. This  implies that
\vspace{-0.2cm}$$\lambda_n(\alpha_{2},\beta_0)<\lambda_n(\alpha_1,\beta_0),\;\; 0\leq n\leq N-1,\;\lambda_n(\alpha_{4},\beta_0)<\lambda_n(\alpha_3,\beta_0),\; 0\leq n\leq N-2,\eqno(3.1)\vspace{-0.2cm}$$
 and in addition,  $\lambda_{N-1}(\alpha_{4},\beta_0)<\lambda_{N-1}(\alpha_3,\beta_0)$ if $\alpha_3\neq\xi$.
Again by (i) of Corollary 4.2 in [22], $\lambda_n(\alpha,\beta_0)$, $0\leq n\leq N-1$,  have the following asymptotic behaviors  near $0$ and $\xi$:
\vspace{-0.2cm}$$\begin{array} {cccc} \lim\limits_{\alpha\rightarrow\pi^-}\lambda_n(\alpha,\beta_0)=\lambda_n(0,\beta_0),\;\;0\leq n \leq N-1,\\[1.0ex]
\lim\limits_{\alpha\rightarrow \xi^-}\lambda_{0}(\alpha,\beta_0)=-\infty,\;\;
\lim\limits_{\alpha\rightarrow \xi^-}\lambda_{n}(\alpha,\beta_0)=\lambda_{n-1}(\xi,\beta_0),\;\;1\leq n\leq N-1,\\[1.0ex]
\lim\limits_{\alpha\rightarrow \xi^+}\lambda_{n}(\alpha,\beta_0)=\lambda_{n}(\xi,\beta_0),\;0\leq n \leq N-2,\;\;\lim\limits_{\alpha\rightarrow\xi^+}\lambda_{N-1}(\alpha,\beta_0)=+\infty.                               \end{array}\vspace{-0.2cm}$$
Thus,
\vspace{-0.2cm}$$\begin{array} {cccc}\lambda_n(\alpha_1,\beta_0)\leq\lambda_n(0,\beta_0)=\lim\limits_{\alpha\rightarrow\pi^-}\lambda_n(\alpha,\beta_0)<\lambda_n(\alpha_4, \beta_0),\; 0\leq n\leq N-1,\\[1.0ex]
\lambda_n(\alpha_3, \beta_0)\leq\lambda_n(\xi, \beta_0)=\lim\limits_{\alpha\rightarrow\xi^-}\lambda_{n+1}(\alpha,\beta_0)<\lambda_{n+1}(\alpha_2,\beta_0),\;0\leq n\leq N-2,\end{array} \vspace{-0.2cm}$$
which together with (3.1) implies that (i) holds. See also Figure 3.1 for $N=4$.

\begin{center}
 \begin{tikzpicture}[scale=0.8]
 \draw [->](-0.5, 0)--(4, 0)node[right]{$\alpha$};
 \draw [->](0, -3.5)--(0, 2.5) node[above]{$\lambda$};
  \draw (0, -1.2).. controls (1,-1.3)  and (1.5,-1.4)..(3.14, -2.2);
 \draw (0, -0.2).. controls (1,-0.3) and (1.5,-0.4)..(3.14, -1.2);
 \draw (0, 1).. controls (1, 0.9) and (1.5, 0.8)..(3.14,  -0.2);
 \draw (0, -2.2).. controls (0.5, -2.25) and (1, -2.3)..(1.4,  -3.5);
  \draw (1.6, 2.5).. controls (2.0, 1.2) and (2.5, 1.05)..(3.14,  1);
   \path(1.5, -3.5) edge [-,dotted] (1.5, 2.5)[line width=0.7pt];
 \path (3.14, -3.5)  edge [-,dotted](3.14, 2.5)[line width=0.7pt];
 \fill (0cm, -2.2cm) circle(2pt);
  \fill (0cm, -1.2cm) circle(2pt);
   \fill (0cm, -0.2cm) circle(2pt);
    \fill (0cm, 1cm) circle(2pt);
     \fill (1.5cm, -1.5cm) circle(2pt);
  \fill (1.5cm, -0.5cm) circle(2pt);
   \fill (1.5cm, 0.65) circle(2pt);
    \node (a) at (0.7,-2.7) {$\lambda_0$};
        \node (a) at (0.7,-1.6) {$\lambda_1$};
        \node (a) at (0.7,-0.6) {$\lambda_2$};
        \node (a) at (0.7,0.6) {$\lambda_3$};
            \node (a) at (2.4,-2.1) {$\lambda_0$};
        \node (a) at (2.4,-1.2) {$\lambda_1$};
        \node (a) at (2.4,0) {$\lambda_2$};
        \node (a) at (2.4,1) {$\lambda_3$};
        \node (a) at (-0.2,-0.2) {$0$};
        \node (a) at (1.5,-0.2) {$\xi$};
          \node (a) at (3.4,-0.2) {$\pi$};
        \draw  (3.14,-2.2) circle(2pt);
            \draw  (3.14,-1.2) circle(2pt);
                \draw  (3.14,-0.2) circle(2pt);
                    \draw  (3.14,1) circle(2pt);
 \end{tikzpicture}
\end{center}
 \begin{center}
   {\small {\bf Figure 3.1.} the $n$-th eigenvalue function of $\alpha$.}
  \end{center}

The proof of assertion (ii) is similar to that for (i) by  (iii) of Corollary 4.2 in [22].

 Now, we show that (iii) holds. Let $\alpha_0\in[0,\xi)\cup(\xi,\pi)$. By (ii) of Corollary 4.2 in [22],
  $\lambda_n(\alpha_0,\beta)$, $0\leq n\leq N-1$,    are strictly increasing in $\beta\in(0,\pi)$  for all $0\leq n \leq N-1$.
Thus,
\vspace{-0.2cm}$$\lambda_n(\alpha_0,\beta_1)<\lambda_n(\alpha_0,\beta_2),\; 0\leq n\leq N-2,\eqno(3.2)\vspace{-0.2cm}$$
and in addition, $\lambda_{N-1}(\alpha_0,\beta_1)<\lambda_{N-1}(\alpha_0,\beta_2)$ if $\beta_2\neq\pi$.
Again by (ii) of Corollary 4.2 in [22],  $\lambda_n(\alpha_0,\beta)$, $0\leq n\leq N-1$, have the following asymptotic behaviors  near $0$ and  $\pi$:
\vspace{-0.2cm}$$\begin{array} {llll}\lim\limits_{\beta\rightarrow \pi^-}\lambda_{n}(\alpha_0,\beta)=\lambda_{n}(\alpha_0,\pi),\;\;0\leq n\leq N-2,\;\;\lim\limits_{\beta\rightarrow \pi^-}\lambda_{N-1}(\alpha_0,\beta)=+\infty,\\[1.0ex]
\lim\limits_{\beta\rightarrow 0^+}\lambda_{0}(\alpha_0,\beta)=-\infty,\;\;\lim\limits_{\beta\rightarrow 0^+}\lambda_{n}(\alpha_0,\beta)=\lambda_{n-1}(\alpha_0,\pi),\;\;1\leq n \leq N-1.
\end{array}\vspace{-0.2cm}$$
Thus,\vspace{-0.2cm}$$\begin{array}{cccc}
\lambda_n(\alpha_0,\beta_2)\leq\lambda_n(\alpha_0,\pi)=\lim\limits_{\beta\rightarrow0^+}\lambda_{n+1}(\alpha_0,\beta)<\lambda_{n+1}(\alpha_0,\beta_1), \;0\leq n\leq N-2,\end{array}\vspace{-0.2cm}$$
which together with (3.2) implies that (iii) holds.
See also Figure 3.2 for $N=4$.

\begin{center}
 \begin{tikzpicture}[scale=0.8]
 \draw [->](-0.5, 0)--(4, 0)node[right]{$\beta$};
 \draw [->](0, -3)--(0, 3) node[above]{$\lambda$};
  \draw (0.1, -3).. controls (0.5,-1.8)  and (2,-1.6)..(3.14, -1.5);
 \draw (0, -1.5).. controls (1,-0.7) and (2,-0.5)..(3.14, -0.2);
 \draw (0, -0.2).. controls (1, 0.5) and (2, 0.7)..(3.14,  1.2);
 \draw (0, 1.2).. controls (1, 1.25) and (2.5, 1.5)..(3.04,  3);
 \path (3.14, -3)  edge [-,dotted](3.14, 3)[line width=0.7pt];
 \fill (3.14cm, -1.5cm) circle(2pt);
  \fill (3.14cm, -0.2cm) circle(2pt);
   \fill (3.14cm, 1.2cm) circle(2pt);
    \node (a) at (1.5,-2.1) {$\lambda_0$};
        \node (a) at (1.5,-1) {$\lambda_1$};
        \node (a) at (1.5,0.3) {$\lambda_2$};
        \node (a) at (1.5,1.2) {$\lambda_3$};
        \node (a) at (-0.2,-0.2) {$0$};
          \node (a) at (3.4,-0.2) {$\pi$};
        \draw  (0,-1.5) circle(2pt);
            \draw  (0,-0.2) circle(2pt);
                \draw  (0,1.2) circle(2pt);
 \end{tikzpicture}.
\end{center}
 \begin{center}
   {\small {\bf Figure 3.2.} the $n$-th eigenvalue function of $\beta$.}
  \end{center}

The proof of assertion (iv)  is similar to that for  (iii) by (iv) of Corollary 4.2 in [22].
This completes the proof.\medskip

The following result is to
compare eigenvalues for an arbitrarily separated BC with those for the   BCs which are Dirichlet at an endpoint.\medskip

\noindent{\bf Corollary 3.1.} {\it Fix a difference equation $\pmb\omega=(1/f,q,w)$ and a separated BC $\mathbf{S}_{\alpha_0,\beta_0}$. Then we have that
\begin{itemize}\vspace{-0.2cm}\item[{\rm (i)}]
 for any $\alpha_0\in(0,\xi)$ and $\beta_0\in(0,\pi)$,
\vspace{-0.2cm}$$\begin{array} {cccc}\lambda_0(\alpha_0,\beta_0)<\{\lambda_0(0,\beta_0),\lambda_0(\alpha_0,\pi)\}<\lambda_1(\alpha_0,\beta_0)
<\{\lambda_1(0,\beta_0),\lambda_1(\alpha_0,\pi)\}<\cdots
 <\\[1.0ex]\lambda_{N-2}(\alpha_0,\beta_0)
 <\{\lambda_{N-2}(0,\beta_0),\lambda_{N-2}(\alpha_0,\pi)\}<\lambda_{N-1}(\alpha_0,\beta_0)<\lambda_{N-1}(0,\beta_0);
 \end{array}\vspace{-0.2cm}$$
\item[{\rm (ii)}] for any $\alpha_0\in(\xi,\pi)$ and $\beta_0\in(0,\pi)$,
\vspace{-0.1cm}$$\begin{array} {cccc}\lambda_0(0,\beta_0)<\lambda_0(\alpha_0,\beta_0)<\{\lambda_1(0,\beta_0),\lambda_0(\alpha_0,\pi)\}<\lambda_1(\alpha_0,\beta_0)
<\{\lambda_2(0,\beta_0),\\[1.0ex]\lambda_1(\alpha_0,\pi)\}<\cdots
 <
 \lambda_{N-2}(\alpha_0,\beta_0)
 <\{\lambda_{N-1}(0,\beta_0),\lambda_{N-2}(\alpha_0,\pi)\}<\lambda_{N-1}(\alpha_0,\beta_0);
 \end{array}\vspace{-0.2cm}$$
\item[{\rm (iii)}]for any $\alpha_0=\xi$ and $\beta_0\in(0,\pi)$,
\vspace{-0.1cm}$$\begin{array} {cccc}\lambda_0(0,\beta_0)<\lambda_0(\alpha_0,\beta_0)<\{\lambda_1(0,\beta_0),\lambda_0(\alpha_0,\pi)\}<\lambda_1(\alpha_0,\beta_0)
<\{\lambda_2(0,\beta_0),\lambda_1(\alpha_0,\pi)\}\\[1.0ex]<\cdots
 <
 \lambda_{N-3}(\alpha_0,\beta_0)
 <\{\lambda_{N-2}(0,\beta_0),\lambda_{N-3}(\alpha_0,\pi)\}<\lambda_{N-2}(\alpha_0,\beta_0)<\lambda_{N-1}(0,\beta_0);
 \end{array}\vspace{-0.2cm}$$
 \item[{\rm (iv)}]for any $\alpha_0\in(0,\xi)$ and $\beta_0=\pi$,
\vspace{-0.2cm}$$\begin{array} {cccc}\lambda_0(\alpha_0,\beta_0)<\lambda_0(0,\beta_0)<\lambda_1(\alpha_0,\beta_0)<\lambda_1(0,\beta_0)
<\cdots<\lambda_{N-2}(\alpha_0,\beta_0)<\lambda_{N-2}(0,\beta_0);
 \end{array}\vspace{-0.2cm}$$
 \item[{\rm (v)}]for any $\alpha_0\in(\xi,\pi)$ and $\beta_0=\pi$,
\vspace{-0.2cm}$$\begin{array} {cccc}\lambda_0(0,\beta_0)<\lambda_0(\alpha_0,\beta_0)<\lambda_1(0,\beta_0)
<\lambda_1(\alpha_0,\beta_0)<\cdots<\lambda_{N-2}(0,\beta_0)<\lambda_{N-2}(\alpha_0,\beta_0);
 \end{array}\vspace{-0.2cm}$$
  \item[{\rm (vi)}]for any $\alpha_0=\xi$ and $\beta_0=\pi$,
\vspace{-0.2cm}$$\begin{array} {cccc}\lambda_0(0,\beta_0)<\lambda_0(\alpha_0,\beta_0)<\cdots<\lambda_{N-3}(0,\beta_0)<\lambda_{N-3}(\alpha_0,\beta_0)<\lambda_{N-2}(0,\beta_0);
 \end{array}\vspace{-0.2cm}$$
  \item[{\rm (vii)}]for any $\alpha_0=0$ and $\beta_0\in(0,\pi)$,
\vspace{-0.2cm}$$\begin{array} {cccc}\lambda_0(\alpha_0,\beta_0)<\lambda_0(\alpha_0,\pi)<\cdots<\lambda_{N-2}(\alpha_0,\beta_0)<\lambda_{N-2}(\alpha_0,\pi)<\lambda_{N-1}(\alpha_0,\beta_0),
 \end{array}\vspace{-0.2cm}$$
\end{itemize}
where the notation $\{\lambda_0(0,\beta_0),\lambda_0(\alpha_0,\pi)\}$ means each of $\lambda_0(0,\beta_0)$ and $\lambda_0(\alpha_0,\pi)$, etc.}\medskip

\noindent{\bf Proof.} (i) and (iii), (i) and (iv), (ii), and (iii) of Theorem 3.1 imply that assertions (i)--(ii), (iii), (iv)--(vi), and (vii) hold, respectively.
 This completes the proof.\medskip

\noindent{\bf 3.2. Inequalities among eigenvalues for different BCs in a natural loop }\medskip

In this subsection, we shall establish  inequalities among eigenvalues for different BCs in a natural loop (given in Lemma 2.7).
We shall  remark that inequalities among eigenvalues for different BCs in a natural loop will play an important role in establishing inequalities
among eigenvalues for coupled BCs and those for some certain separated ones, and among eigenvalues for different coupled BCs in subsections 3.3 and 3.4.

Firstly, we shall establish  inequalities among eigenvalues for different BCs in the natural loops $\mathcal{C}_{1,4,z, b_{21}}$ and
$\mathcal{C}_{1,4,z, a_{12}}$, separately.\medskip

\noindent{\bf Theorem 3.2.} {\it
 Fix a difference equation $\pmb\omega=(1/f,q,w)$. Let
 \vspace{-0.1cm}$$\mathbf{A}(a_{12},b_{21}):=\left [\begin{array} {cccc}1&a_{12}&\bar{z}&0\\
0&z&b_{21}&1\end{array}  \right ]\in\mathcal{O}_{1,4}^{\mathbb{C}}.\vspace{-0.1cm}$$
 Then $(\pmb\omega, \mathbf{A}(a_{12},b_{21}))$  has exactly $N$ eigenvalues if  $a_{12}\neq1/f_0$ and exactly $N-1$ eigenvalues if  $a_{12}=1/f_0$;  $(\pmb\omega, \mathbf{S_1})$
  has exactly $N$ eigenvalues   in any case;
  $(\pmb\omega, \mathbf{S_2})$
  has exactly $N-1$ eigenvalues if  $a_{12}\neq1/f_0$ and exactly $N-2$ eigenvalues   if  $a_{12}=1/f_0$,
 where $\mathbf{S}_1 $ and $\mathbf{S}_2$ are specified in Lemma {\rm2.6}. Further, for any  $a_{12}^{(1)}<a_{12}^{(2)}\leq 1/f_0<a_{12}^{(3)}<a_{12}^{(4)}$ and $b_{21}^{(1)}<b_{21}^{(2)}$, we have that
  \begin{itemize}\vspace{-0.2cm}
\item[{\rm (i)}] the eigenvalues $\lambda_n(a_{12}^{(i)})$ of the SLPs $(\pmb\omega, \mathbf{A}(a_{12}^{(i)},b_{21}))$, $i=1,\cdots, 4$, and $\lambda_n(\mathbf{S}_1)$ of $(\pmb\omega, \mathbf{S}_1)$
satisfy the following inequalities:
 \vspace{-0.1cm}$$\vspace{-0.2cm}\begin{array} {cccc}\lambda_0(a_{12}^{(3)})\leq\lambda_0(a_{12}^{(4)})\leq\lambda_0(\mathbf{S}_1)\leq\lambda_0(a_{12}^{(1)})\leq\lambda_0(a_{12}^{(2)})\leq
\lambda_1(a_{12}^{(3)})\leq\lambda_1(a_{12}^{(4)})\leq
\lambda_1(\mathbf{S}_1)\\[1.0ex]\leq\lambda_1(a_{12}^{(1)})\leq\lambda_1(a_{12}^{(2)})\leq
\cdots\leq
\lambda_{N-2}(a_{12}^{(3)})\leq\lambda_{N-2}(a_{12}^{(4)})\leq\lambda_{N-2}(\mathbf{S}_1)\leq\lambda_{N-2}(a_{12}^{(1)})\\[1.0ex]\leq\lambda_{N-2}(a_{12}^{(2)})\leq
\lambda_{N-1}(a_{12}^{(3)})\leq\lambda_{N-1}(a_{12}^{(4)})\leq\lambda_{N-1}(\mathbf{S}_1)\leq\lambda_{N-1}(a_{12}^{(1)}),\end{array}\vspace{-0.01cm}$$
and  in addition, $\lambda_{N-1}(a_{12}^{(1)})\leq\lambda_{N-1}(a_{12}^{(2)})$ if $a_{12}^{(2)}<1/f_0$;
\item[{\rm (ii)}] the eigenvalues $\lambda_n(b_{21}^{(j)})$ of the SLPs $(\pmb\omega, \mathbf{A}(a_{12},b_{21}^{(j)}))$, $j=1,2$, and $\lambda_n(\mathbf{S}_2)$ of $(\pmb\omega, \mathbf{S}_2)$
satisfy the following inequalities:
\vspace{-0.1cm}$$\begin{array} {cccc}\vspace{-0.2cm}\lambda_0(b_{21}^{(1)})\leq\lambda_0(b_{21}^{(2)})\leq\lambda_0(\mathbf{S}_2)\leq\lambda_1(b_{21}^{(1)})\leq
\lambda_1(b_{21}^{(2)})\leq\lambda_1(\mathbf{S}_2)\leq\cdots\leq\\[2ex]
\lambda_{N-3}(b_{21}^{(1)})\leq\lambda_{N-3}(b_{21}^{(2)})\leq\lambda_{N-3}(\mathbf{S}_2)\leq
\lambda_{N-2}(b_{21}^{(1)})\leq\lambda_{N-2}(b_{21}^{(2)}),\end{array}\vspace{-0.15cm}$$
and in addition, $\lambda_{N-2}(b_{21}^{(2)})\leq
\lambda_{N-2}(\mathbf{S}_2)\leq\lambda_{N-1}(b_{21}^{(1)})\leq\lambda_{N-1}(b_{21}^{(2)})$ if $a_{12}\neq1/f_0$.
    \end{itemize}}\medskip

\noindent{\bf Proof.} The number of eigenvalues of $(\pmb\omega,\mathbf{A}(a_{12},b_{21}))$, $(\pmb\omega,\mathbf{S}_1)$, and $(\pmb\omega,\mathbf{S}_2)$ can be obtained by Lemma 2.4 and direct computations.

Let
$\mathbf{A}(s)$ and $\mathcal{C}_{1,4,z,b_{21}}=\{\mathbf{A}(s):s\in\mathbb{R}\}\cup\{\mathbf{S}_1\}$ be given as that in (i) of Lemma 2.7.
Then  $\mathbf{A}(a_{12}^{(i)})=\mathbf{A}(a_{12}^{(i)},b_{21})$, $i=1,\cdots,4$. By (i)--(ii) of Theorem 4.1 in [22],  the eigenvalue functions $\lambda_n(\mathbf{A}(s))$ are continuous and
non-decreasing in $(-\infty,1/ f_0)$ and $(1/ f_0,+\infty)$ for all $0\leq n\leq  N-1$. Thus,  one gets that
 \vspace{-0.2cm}$$\lambda_n(a_{12}^{(1)})\leq\lambda_n(a_{12}^{(2)}),\;0\leq n\leq N-2, \;\lambda_n(a_{12}^{(3)})\leq\lambda_n(a_{12}^{(4)}),\;0\leq n\leq N-1, \eqno(3.3)\vspace{-0.2cm}$$
and in addition, $\lambda_{N-1}(a_{12}^{(1)})\leq\lambda_{N-1}(a_{12}^{(2)})$ if $a_{12}^{(2)}<1/f_0$.    By (iii) of Theorem 4.1 in [22],
$\lambda_n(\mathbf{A}(s))$, $0\leq n\leq N-1$,  have  asymptotic behaviors     near $1/ f_0$ as follows:
$$\vspace{-0.1cm}\begin{array} {llll}\lim\limits_{s\rightarrow(1/ f_0)^-}\lambda_n({\mathbf{A}}(s))=\lambda_{n}({\mathbf{A}}(1/ f_0)),\;0\leq n \leq N-2,\;\lim\limits_{s\rightarrow(1/ f_0)^-}\lambda_{N-1}({\mathbf{A}}(s))=+\infty,
\\[1.0ex]
\lim\limits_{s\rightarrow(1/ f_0)^+}\lambda_{0}({\mathbf{A}}(s))=-\infty,
\lim\limits_{s\rightarrow(1/ f_0)^+}\lambda_n({\mathbf{A}}(s))=\lambda_{n-1}({\mathbf{A}}(1/ f_0)),
  \;1\leq n \leq N-1.
\end{array}\vspace{0.05cm}$$
Since  $\mathcal{C}_{1,4,z,b_{21}}\backslash\{\mathbf{A}(1/f_0)\}$ is connected by Remark 2.1 and $(\pmb\omega,\mathbf{A})$ has
 exactly $N$ eigenvalues for each $\mathbf{A}\in\mathcal{C}_{1,4,z,b_{21}}\backslash\{\mathbf{A}(1/f_0)\}$,  $\lambda_n$ restricted in $\mathcal{C}_{1,4,z,b_{21}}\backslash\{\mathbf{A}(1/f_0)\}$ is continuous for each $0\leq n\leq N-1$ by Lemma 2.5.
This, together with Lemma 2.6, implies that
$\lim\limits_{s\rightarrow \pm\infty}\lambda_{n}({\mathbf{A}}(s))=\lambda_{n}(\mathbf{S}_1)$ for all $0\leq n \leq N-1$.
See Figure 3.3 for $N=4$.
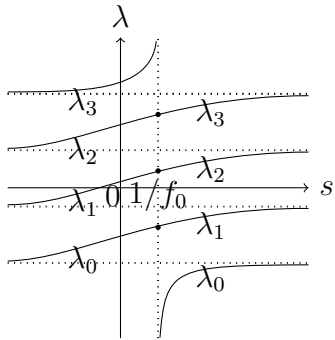
\begin{center}
 \begin{tikzpicture}[scale=0.5]
 \draw [->](-3, 0)--(5, 0)node[right]{$s$};
 \draw [->](0, -4)--(0, 4) node[above]{$\lambda$};
  \draw (-3, 1.05).. controls (-1,1.06)  and (1,2.44)..(5, 2.45);
 \draw (-3, -0.45).. controls (-1,-0.48) and (1,0.94)..(5, 0.95);
 \draw (-3, -1.95).. controls (-1, -1.89) and (1, -0.52)..(5,  -0.55);
 \draw (-3, 2.55).. controls (-1, 2.56) and (0.8, 2.53)..(0.95,  3.9);
  \draw (1.05, -3.99).. controls (1.2,-2.1) and (1.4, -2.07)..(5,  -2.05);
   \path(1, -4) edge [-,dotted] (1, 4)[line width=0.7pt];
 \path (-3, 2.5)  edge [-,dotted](5, 2.5)[line width=0.7pt];
  \path (-3, 1)  edge [-,dotted](5, 1)[line width=0.7pt];
   \path (-3, 2.5)  edge [-,dotted](5, 2.5)[line width=0.7pt];
      \path (-3, -0.5)  edge [-,dotted](5, -0.5)[line width=0.7pt];
   \path (-3, -2)  edge [-,dotted](5, -2)[line width=0.7pt];
 \fill (1cm, 1.95cm) circle(2pt);
  \fill (1cm, 0.45cm) circle(2pt);
   \fill (1cm, -1.05cm) circle(2pt);
    \node (a) at (-1,-1.9) {$\lambda_0$};
        \node (a) at (-1,-0.4) {$\lambda_1$};
        \node (a) at (-1,1) {$\lambda_2$};
        \node (a) at (-1,2.3) {$\lambda_3$};
            \node (a) at (2.4,-2.4) {$\lambda_0$};
        \node (a) at (2.4,-1) {$\lambda_1$};
        \node (a) at (2.4,0.5) {$\lambda_2$};
        \node (a) at (2.4,2) {$\lambda_3$};
        \node (a) at (-0.2,-0.2) {$0$};
        \node (a) at (1,-0.2) {$1/f_0$};
 \end{tikzpicture}
\end{center}\vspace{-0.2cm}
 \begin{center}\vspace{-0.2cm}
   {\small {\bf Figure 3.3.} the $n$-th eigenvalue function of $s$.}
  \end{center}\vspace{-0.2cm}
Thus,
\vspace{-0.2cm}$$\begin{array}{cccc}\lambda_n(a_{12}^{(4)})\leq\lambda_n(\mathbf{S}_1)\leq\lambda_n(a_{12}^{(1)}),\;0\leq n\leq N-1,\\[1.0ex]
\lambda_n(a_{12}^{(2)})\leq\lambda_n(\mathbf{A}(1/f_0))=\lim\limits_{s\rightarrow(1/f_0)^+}\lambda_{n+1}(\mathbf{A}(s))
\leq\lambda_{n+1}(a_{12}^{(3)}),\;0\leq n\leq N-2.\end{array}\eqno(3.4)\vspace{-0.2cm}$$
 Hence (3.3)--(3.4) implies (i) holds.

Then we show that  (ii) holds.
Let
$\hat{\mathbf{A}}(t)$ and $\mathcal{C}_{1,4,z,a_{12}}=\{\hat{\mathbf{A}}(t):t\in\mathbb{R}\}\cup\{\mathbf{S}_2\}$ be given as that in (i) of Lemma 2.7.
Then  $\mathbf{A}(a_{12},b_{21}^{(j)})=\hat{\mathbf{A}}(b_{21}^{(j)})$, $j=1,2$.

Let $a_{12}\neq1/f_0$. Then $\lambda_n(\hat{\mathbf{A}}(t))$ are continuous and non-decreasing in $t\in\mathbb{R}$ for all $0\leq n\leq N-1$   by (i)--(ii) of Theorem 4.1 in [22]. Thus,  for each $ 0\leq n\leq N-1$, one has that
  \vspace{-0.2cm}$$\lambda_n(b_{21}^{(1)})\leq\lambda_n(b_{21}^{(2)}).  \eqno(3.5)\vspace{-0.2cm}$$

 To see the limits of $\lambda_n(\hat{\mathbf{A}}(t))$ at $\pm\infty$, we notice that for $t\neq 0$,
$$\vspace{-0.2cm}\hat{\mathbf{A}}(t)=\left [\begin{array} {cccc}1&a_{12}&\bar{z}&0\\
0&z&t&1\end{array}  \right ]=\left [\begin{array} {cccc}1&a_{12}-z\bar{z}/t&0&-\bar{z}/t\\
0&-z/t&-1&-1/t\end{array}  \right ].\eqno(3.6)\vspace{-0.2cm}$$

In the case that $a_{12}>1/f_0$, direct computations show that   $\hat{\mathbf{A}}(t)\in\mathcal{B}_{1,3r}^+$ if $t<0$; $\hat{\mathbf{A}}(t)\in\mathcal{B}_{1,3}^-$ if $t>0$; and $\mathbf{S}_2\in\mathcal{B}_{1,3r}$, where
$$\vspace{-0.2cm}\begin{array} {cccc}\mathcal{B}_{1,3r}^+:=\left\{ \mathbf{A}\in\mathcal{O}_{1,3}^{\mathbb{C}}: a_{12}\geq 1/ f_0,b_{22}\geq 0, (a_{12}-1/f_0)b_{22}>|z|^2\right\},\\[1.0ex]\mathcal{B}_{1,3}^-:=\left\{ \mathbf{A}\in\mathcal{O}_{1,3}^{\mathbb{C}}:  (a_{12}-1/ f_0)b_{22}<|z|^2\right\},\; \mathbf{C}:=\left [\begin{array} {cccc}1&1/ f_0&0&0\\
0&0&1&0\end{array}  \right ],\\[2.0ex]
\mathcal{B}_{1,3r}:=\left\{\mathbf{A}\in\mathcal{O}_{1,3}^{\mathbb{C}}:  (a_{12}-1/ f_0)b_{22}=|z|^2, a_{12}\geq 1/ f_0,b_{22}\geq 0\right\}\backslash\left\{\mathbf{C}\right\}.\;\end{array}\vspace{-0.05cm}$$
Note that $\lim\limits_{t\rightarrow\pm\infty}\hat{\mathbf{A}}(t)=\mathbf{S}_2$ by Lemma 2.6. Then, it follows from (iiia) of Theorem 4.3 in [22] that
\vspace{-0.2cm}$$\vspace{-0.2cm}\begin{array} {llll}\lim\limits_{t\rightarrow-\infty}\lambda_0(\hat{\mathbf{A}}(t))=-\infty,\;
\lim\limits_{t\rightarrow-\infty}\lambda_n(\hat{\mathbf{A}}(t))=\lambda_{n-1}({\mathbf{S}_2}),\;1\leq n \leq N-1,\\[1.0ex]
\lim\limits_{t\rightarrow+\infty}\lambda_n(\hat{\mathbf{A}}(t))=\lambda_{n}({\mathbf{S}_2}),\;0\leq n \leq N-2,\;
\lim\limits_{t\rightarrow+\infty}\lambda_{N-1}(\hat{\mathbf{A}}(t))=+\infty.
\end{array}$$
See Figure 3.4 for $N=4$.
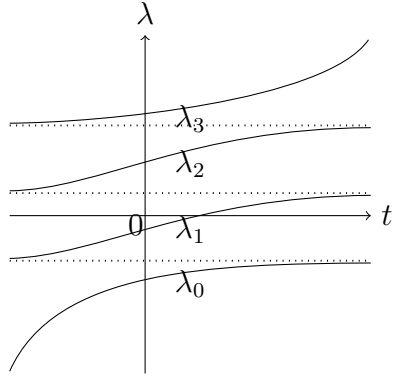
\begin{center}
 \begin{tikzpicture}[scale=0.6]
 \draw [->](-3, 0)--(5, 0)node[right]{$t$};
 \draw [->](0, -3.5)--(0, 4) node[above]{$\lambda$};
 \draw (-3, 0.55).. controls (-1,0.58) and (1, 1.94)..(5, 1.95);
 \draw (-3, -0.95).. controls (-1, -0.89) and (1, 0.42)..(5,  0.45);
 \draw (-3, 2.05).. controls (-1, 2.06) and (4, 2.53)..(4.95,  3.9);
  \draw (-3, -3.45).. controls (-2,-1.2) and (1.4, -1.07)..(5,  -1.05);
 \path (-3, 2)  edge [-,dotted](5, 2) [line width=0.7pt];
  \path (-3, 0.5)  edge [-,dotted](5, 0.5) [line width=0.7pt];
      \path (-3, -1)  edge [-,dotted](5, -1) [line width=0.7pt];
    \node (a) at (1,-1.5) {$\lambda_0$};
        \node (a) at (1,-0.3) {$\lambda_1$};
        \node (a) at (1,1.15) {$\lambda_2$};
        \node (a) at (1,2.15) {$\lambda_3$};
        \node (a) at (-0.2,-0.2) {$0$};
 \end{tikzpicture}
\end{center}\vspace{-0.2cm}
 \begin{center}\vspace{-0.2cm}
   {\small {\bf Figure 3.4.} the $n$-th eigenvalue function of $t$.}
  \end{center}\vspace{-0.2cm}
Thus,
\vspace{-0.2cm}$$\lambda_n(b_{21}^{(2)})\leq\lim\limits_{t\rightarrow+\infty}\lambda_n(\hat{\mathbf{A}}(t))=
\lambda_n(\mathbf{S}_2)=
\lim\limits_{t\rightarrow-\infty}\lambda_{n+1}(\hat{\mathbf{A}}(t))\leq
\lambda_{n+1}(b_{21}^{(1)}),\eqno(3.7)\vspace{-0.2cm}$$
for  $0\leq n\leq N-2$. (3.5) and (3.7) implies (ii) holds in the case that $a_{12}>1/f_0$.

In the case that $a_{12}<1/f_0$, direct computations imply that $\hat{\mathbf{A}}(t)\in\mathcal{B}_{1,3}^-$ if $t<0$;  $\hat{\mathbf{A}}(t)\in\mathcal{B}_{1,3l}^+$ if $t>0$; and $\mathbf{S}_2\in\mathcal{B}_{1,3l}$, where
\vspace{-0.2cm}$$\begin{array}{cccc}\mathcal{B}_{1,3l}^+:=\left\{ \mathbf{A}\in\mathcal{O}_{1,3}^{\mathbb{C}}:  a_{12}\leq 1/ f_0,b_{22}\leq 0,  (a_{12}-1/ f_0)b_{22}>|z|^2\right\},\\[1.0ex]\mathcal{B}_{1,3l}:=\left\{ \mathbf{A}\in\mathcal{O}_{1,3}^{\mathbb{C}}:  (a_{12}-1/ f_0)b_{22}=|z|^2, a_{12}\leq 1/ f_0,b_{22}\leq 0\right\}\backslash\left\{\mathbf{C}\right\}.\end{array}\vspace{-0.2cm}$$
By (iiib) of Theorem 4.3 in [22], similar arguments above yield that (ii) holds in this case.

Let $a_{12}=1/f_0$ Then $\mathbf{S}_2=\mathbf{C}$. Since $\mathcal{C}_{1,4,z,a_{12}}\backslash\{\mathbf{S}_2\}$ is connected by Remark 2.1 and  $(\pmb\omega,\mathbf{A})$ has  exactly $N-1$ eigenvalues for each $\mathbf{A}\in\mathcal{C}_{1,4,z,a_{12}}\backslash\{\mathbf{S}_2\}$, by Lemma 2.5 the eigenvalue function $\lambda_n$ is continuous and locally
forms a continuous eigenvalue branch in $\mathcal{C}_{1,4,z, a_{12}}\backslash\{\mathbf{S}_2\}$ for each $0\leq n\leq N-2$. By Theorem 4.6 in [23], $\lambda_n(\hat{\mathbf{A}}(t))$ is non-decreasing in $t\in\mathbb{R}$, and  thus (3.5) holds for each $0\leq n\leq N-2$.
 By (3.6) and direct computations, it follows  that $\hat{\mathbf{A}}(t)\in\mathcal{B}_{1,3r}$ if $t<0$; and $\hat{\mathbf{A}}(t)\in\mathcal{B}_{1,3l}$ if $t>0$.
In addition,  $\lim\limits_{t\rightarrow\pm\infty}\hat{\mathbf{A}}(t)=\mathbf{S}_2$ by Lemma 2.6. It follows from (iiic) of Theorem 4.3 in [22] that
\vspace{-0.2cm}$$\begin{array} {llll}\lim\limits_{t\rightarrow-\infty}\lambda_0(\hat{\mathbf{A}}(t))=-\infty,
\lim\limits_{t\rightarrow-\infty}\lambda_n(\hat{\mathbf{A}}(t))=\lambda_{n-1}({\mathbf{S}_2}),\;1\leq n \leq N-2,\\[1.0ex]
\lim\limits_{t\rightarrow+\infty}\lambda_n(\hat{\mathbf{A}}(t))=\lambda_{n}({\mathbf{S}_2}),\;0\leq n \leq N-3,
\lim\limits_{t\rightarrow+\infty}\lambda_{N-2}(\hat{\mathbf{A}}(t))=+\infty.
\end{array}\vspace{-0.1cm}$$
Thus, (3.7) holds for each $0\leq n\leq N-3$.
 Hence, (ii) holds.  The proof is  complete.\medskip

\noindent{\bf Remark 3.1.} The inequalities in Theorem 3.2 may not be strict. See the following example. \medskip

\noindent{\bf Example 3.1.}  Consider (1.1)--(1.2), where
\vspace{-0.2cm}$$
f_0 =1, \ f_1 =1,\ f_2=1, \; q_1 =q_2 =0, \; w_1 =w_2 =1,\; N=2,
\vspace{-0.2cm}$$
and
\vspace{-0.2cm}$$
\mathbf A_1(s): =\left[\begin{array}{cccc} 1 & s &     -1 & 0 \\
                                 0 &     -1 & 0 & 1 \end{array}\right]
\in \mathcal O_{1,4}^{\mathbb C}.
\vspace{-0.2cm}$$
Then, by Lemma 2.3,
\vspace{-0.2cm}$$
\Gamma(\lambda) =-(s-1) \lambda^2
+2(s-2) \lambda.
\vspace{-0.05cm}$$
Thus, for each $s\in(-\infty,1)\cup(1,\infty)$, there are exactly  two eigenvalues for $\mathbf A_1(s)$ and exactly one eigenvalue for $\mathbf A_1(1)$:
\vspace{-0.2cm}$$\begin{array}{ll}\vspace{0.2cm}
\lambda_0(s)  =
\begin{cases}                      0 & \text{ if } s\leq1, \\
       2(s-2)/(s-1) & \text{ if } 1<s\leq2,\\
       0 & \text{ if } s>2,\end{cases}\;\;
\lambda_1(s)  =
\begin{cases} 2(s-2)/(s-1) & \text{ if } s<1, \\
                            0 & \text{ if } 1<s\leq2,\\
       2(s-2)/(s-1) & \text{ if } s>2.\end{cases}
\end{array}\vspace{-0.05cm}$$
Note that \vspace{-0.05cm}$$\mathbf{S}_1=\lim_{s\rightarrow+\infty}\mathbf{A}_1(s)=\left [\begin{array} {cccc}0&1&0&0\\
0&0&0&1\end{array}  \right ].$$
It is easy to see  that there are exactly two eigenvalues for $\mathbf{S}_1$, and they are $0$ and $2$. See Figure 3.5.
\begin{center}
 \begin{tikzpicture}[scale=0.3,domain=-10:0.75,smooth,variable=\x]
 \draw [->](-10, 0)--(11.5, 0)node[right]{$s$};
 \draw [->](0, -8)--(0, 10) node[above]{$\lambda$};
\draw plot(\x, {2*(\x-2)/(\x-1)}) [line width=1.5pt];
\draw plot(\x+11.2, {2*(\x+11.2-2)/(\x+11.2-1)})[line width=1.5pt];
\draw plot (\x,0)[line width=1.5pt];
\draw plot (\x+10.7,0)[line width=1.5pt];
 \path (1, -8)  edge [-,dotted](1, 10)[line width=1pt];
    \node (a) at (5,-1) {$\lambda_0$};
        \node (a) at (2.2,-5) {$\lambda_0$};
        \node (a) at (-5,-1) {$\lambda_0$};
        \node (a) at (7,2.5) {$\lambda_1$};
        \node (a) at (1.7,1) {$\lambda_1$};
        \node (a) at (-7,3) {$\lambda_1$};
            \node (a) at (-0.5,-0.7) {$0$};
              \node (a) at (1,-0.7) {$1$};
                \node (a) at (2.3,-0.7) {$2$};
                 \node (a) at (0,2) {$-$};
                  \node (a) at (-0.7,2) {$2$};
 \end{tikzpicture}
\end{center}\vspace{-0.2cm}
 \begin{center}\vspace{-0.2cm}
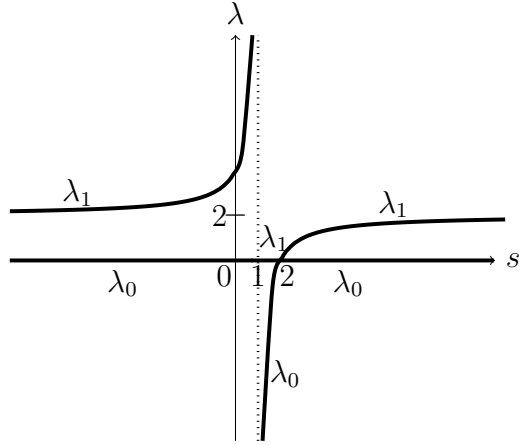
   {\small {\bf Figure 3.5.} the $n$-th eigenvalue function of $s$ in Example 3.1.}
  \end{center}\vspace{-0.2cm}
Then  $\lambda_0(\mathbf{S}_1)=\lambda_0(\mathbf{A}_1(s))<\lambda_1(\mathbf{S}_1)<\lambda_1(\mathbf{A}_1(s))$ for
$s<1$;
 $\lambda_0(\mathbf{S}_1)=\lambda_0(\mathbf{A}_1(s))<\lambda_1(\mathbf{S}_1)$ for $s=1$;
 $\lambda_0(\mathbf{A}_1(s))<\lambda_0(\mathbf{S}_1)=\lambda_1(\mathbf{A}_1(s))<\lambda_1(\mathbf{S}_1)$ for $1< s<2$;
 $\lambda_0(\mathbf{A}_1(s))=\lambda_0(\mathbf{S}_1)=\lambda_1(\mathbf{A}_1(s))<\lambda_1(\mathbf{S}_1)$ for $s=2$; and
$\lambda_0(\mathbf{A}_1(s))=\lambda_0(\mathbf{S}_1)<\lambda_1(\mathbf{A}_1(s))<\lambda_1(\mathbf{S}_1)$ for $s>2$.
\medskip

Secondly, we shall establish  inequalities among eigenvalues for different BCs  in the natural loops  $\mathcal{C}_{2,4,z, b_{21}}$ and
$\mathcal{C}_{2,4,z, a_{11}}$, separately.\medskip

\noindent{\bf Theorem 3.3.} {\it  Fix a difference equation $\pmb\omega=(1/f,q,w)$. Let
\vspace{-0.2cm}$$\mathbf{A}(a_{11},b_{21}):=\left [\begin{array} {cccc}a_{11}&-1&\bar{z}&0\\
z&0&b_{21}&1\end{array}  \right ] \in \mathcal{O}_{2,4}^{\mathbb{C}}.\vspace{-0.2cm}$$
Then similar results in Theorem {\rm 3.2} hold with $a_{12}$,  $a_{12}^{(i)}$,   $1/ f_0$, $\mathbf{S}_1$, and
$\mathbf{S}_2$ replaced by $a_{11}$, $a_{11}^{(i)}$, $- f_0$, $\mathbf{S}_3$, and $\mathbf{S}_4$, separately, where $i=1,\cdots,4$, and $\mathbf{S}_3$ and $\mathbf{S}_4$ are specified in Lemma {\rm2.6}.}\medskip

\noindent{\bf Proof.} By a similar method  to that used in the proof of Theorem 3.2,  one can show that Theorem 3.3 holds with the help of Theorems 4.2 and 4.4 in [22]. \medskip

Thirdly, we shall establish  inequalities among eigenvalues for different BCs  in the natural loops  $\mathcal{C}_{2,3,z, b_{22}}$ and
$\mathcal{C}_{2,3,z, a_{11}}$, separately. We shall remark that here we only give the inequalities  in the case that $z\neq0$ since we shall apply Theorem 3.4 to coupled BCs, which satisfy that $z\neq0$. One can establish the inequalities  in the case that $z=0$ with a similar  method.\medskip

\noindent{\bf Theorem 3.4.} {\it Fix a difference equation $\pmb\omega=(1/f,q,w)$. Let
\vspace{-0.2cm}$$\mathbf{A}(a_{11},b_{22}):=\left [\begin{array} {cccc}a_{11}&-1&0&\bar{z}\\
z&0&-1&b_{22}\end{array}  \right ]\in \mathcal{O}_{2,3}^{\mathbb{C}},\vspace{-0.2cm}$$
where $z\neq0$. Then $(\pmb\omega,\mathbf{A}(a_{11},b_{22}))$ has
exactly $N$ eigenvalues  if $ b_{22}(a_{11}+f_0)\neq|z|^2$, and exactly $N-1$ eigenvalues  if $ b_{22}(a_{11}+f_0)=|z|^2$;
 $(\pmb\omega,\mathbf{S}_5)$ has  exactly $N$ eigenvalues  if $ b_{22}\neq0$, and exactly $N-1$ eigenvalues  if $ b_{22}=0$;
  $(\pmb\omega,\mathbf{S}_6)$ has exactly $N$ eigenvalues  if $ a_{11}+f_0\neq0$, and exactly $N-1$ eigenvalues  if $ a_{11}+f_0=0$, where $\mathbf{S}_5 $ and $\mathbf{S}_6$ are specified in Lemma {\rm2.6}. Further, we have that \begin{itemize}\vspace{-0.2cm}
\item[{\rm (i)}]  in the case that  $b_{22}=0$, for any $a_{11}^{(1)}< a_{11}^{(2)}$,  the eigenvalues $\lambda_n(a_{11}^{(i)})$ of $(\pmb\omega, \mathbf{A}(a_{11}^{(i)},b_{22}))$, $i=1,2$, and $\lambda_n(\mathbf{S}_5)$ of  $(\pmb\omega,\mathbf{S}_5)$ satisfy the following inequalities:
\vspace{-0.2cm}$$\begin{array} {cccc}\lambda_0(a_{11}^{(1)})\leq\lambda_0(a_{11}^{(2)})\leq\lambda_0(\mathbf{S}_5)\leq\lambda_1(a_{11}^{(1)})\leq\lambda_1(a_{11}^{(2)})
\leq\lambda_1(\mathbf{S}_5)\leq\cdots\leq\\[1.0ex]
\lambda_{N-2}(a_{11}^{(1)})\leq\lambda_{N-2}(a_{11}^{(2)})\leq\lambda_{N-2}(\mathbf{S}_5)\leq\lambda_{N-1}(a_{11}^{(1)})\leq\lambda_{N-1}(a_{11}^{(2)});
\end{array}\eqno(3.8)\vspace{-0.2cm}$$
\item[{\rm (ii)}] in the case that   $b_{22}\neq0$,  for any $a_{11}^{(1)}<a_{11}^{(2)}\leq|z|^2/b_{22}-f_0$ and $|z|^2/b_{22}-f_0<a_{11}^{(3)}<a_{11}^{(4)}$, the eigenvalues $\lambda_n(a_{11}^{(i)})$ of $(\pmb\omega, \mathbf{A}(a_{11}^{(i)},b_{22}))$, $i=1,\cdots, 4$, and $\lambda_n(\mathbf{S}_5)$ of $(\pmb\omega,\mathbf{S}_5)$ satisfy the following inequalities:
\vspace{-0.2cm}$$\vspace{-0.2cm}\begin{array} {cccc} \lambda_0(a_{11}^{(3)})\leq \lambda_0(a_{11}^{(4)})\leq\lambda_0(\mathbf{S}_5)\leq\lambda_0(a_{11}^{(1)})\leq\lambda_0(a_{11}^{(2)})\leq\\[1.0ex]
\lambda_1(a_{11}^{(3)})\leq\lambda_1(a_{11}^{(4)})
\leq\lambda_1(\mathbf{S}_5)\leq\lambda_1(a_{11}^{(1)})\leq \lambda_1(a_{11}^{(2)})\\[1.0ex]
\leq\cdots\leq\lambda_{N-2}(a_{11}^{(3)})\leq \lambda_{N-2}(a_{11}^{(4)})\leq\lambda_{N-2}(\mathbf{S}_5)\leq\lambda_{N-2}(a_{11}^{(1)})\leq\\[1.0ex]\lambda_{N-2}(a_{11}^{(2)})\leq
\lambda_{N-1}(a_{11}^{(3)})\leq\lambda_{N-1}(a_{11}^{(4)})\leq\lambda_{N-1}(\mathbf{S}_5)\leq\lambda_{N-1}(a_{11}^{(1)}),\end{array}\eqno(3.9)\vspace{-0.2cm}$$\vspace{-0.1cm}
and in addition,  $\lambda_{N-1}(a_{11}^{(1)})\leq\lambda_{N-1}(a_{11}^{(2)})$ if $a_{11}^{(2)}<|z|^2/b_{22}-f_0$;
 \item[{\rm (iii)}] in the case that  $a_{11}+f_0=0$, for any $b_{22}^{(1)}< b_{22}^{(2)}$,  the eigenvalues $\lambda_n(b_{22}^{(i)})$ of $(\pmb\omega, \mathbf{A}(a_{11},b_{22}^{(i)}))$ and  $\lambda_n(\mathbf{S}_6)$ of $(\pmb\omega,\mathbf{S}_6)$ satisfy {\rm(3.8)} with $a_{11}^{(i)}$ and $\mathbf{S}_5$ replaced by $b_{22}^{(i)}$ and  $\mathbf{S}_6$, separately, where  $i=1,2$;
\item[{\rm (iv)}] in the case that  $a_{11}+f_0\neq0$, for any $b_{22}^{(1)}<b_{22}^{(2)}\leq|z|^2/(a_{11}+f_0)$ and $|z|^2/(a_{11}+f_0)<b_{22}^{(3)}<b_{22}^{(4)}$,  the eigenvalues $\lambda_n(b_{22}^{(i)})$ of $(\pmb\omega, \mathbf{A}(a_{11},b_{22}^{(i)}))$  and  $\lambda_n(\mathbf{S}_6)$ of $(\pmb\omega,\mathbf{S}_6)$ satisfy {\rm(3.9)} with $a_{11}^{(i)}$ and $\mathbf{S}_5$ replaced by $b_{22}^{(i)}$ and $\mathbf{S}_6$, separately, where $i=1,\cdots, 4$,
and in addition, $\lambda_{N-1}(b_{22}^{(1)})\leq\lambda_{N-1}(b_{22}^{(2)})$ if $b_{22}^{(2)}<|z|^2/(a_{11}+f_0)$.
     \end{itemize}}\medskip

\noindent{\bf Proof.} By a similar method to that used in the proof of (ii) in Theorem 3.2,  one can show that (i)  holds  with the help of Theorems 4.3--4.4 of [22];    (iii) holds with the help of Theorems 4.2 and 4.4 of [22].
 By a similar method to that  used in the proof of (i) in Theorem 3.2, one gets that (ii) and (iv) hold with the help of Theorem 4.4 of [22]. This completes the proof.\medskip

Fourthly, we shall establish  inequalities among eigenvalues for different BCs  in the natural loops $\mathcal{C}_{1,3,z, b_{22}}$ and
$\mathcal{C}_{1,3,z, a_{12}}$ with $z\neq0$, separately.\medskip

\noindent{\bf Theorem 3.5.} {\it  Fix a difference equation $\pmb\omega=(1/f,q,w)$. Let
\vspace{-0.2cm}$$\mathbf{A}(a_{12},b_{22}):=\left [\begin{array} {cccc}1&a_{12}&0&\bar{z}\\
0&z&-1&b_{22}\end{array}  \right ]\in \mathcal{O}_{1,3}^{\mathbb{C}},\vspace{-0.2cm}$$
where $z\neq0$. Then similar results in Theorem {\rm 3.4} hold for $a_{11}$, $a_{11}^{(i)}$, $a_{11}+f_0$, $|z|^2/b_{22}-f_0$, $\mathbf{S}_5$, and $\mathbf{S}_6$ replaced by $a_{12}$,  $a_{12}^{(i)}$, $a_{12}-1/f_0$, $|z|^2/b_{22}+1/f_0$, $\mathbf{S}_7$, and  $\mathbf{S}_8$, separately, where $i=1,\cdots, 4$, and $\mathbf{S}_7$ and $\mathbf{S}_8$ are specified in Lemma {\rm2.6}.}\medskip

\noindent{\bf Proof.} By a similar method to that used in the proof of (ii) in Theorem 3.2,  one can show that (i)  holds  with the help of Theorems 4.3--4.4 of [22], and (iii) holds with the help of Theorems 4.1 and 4.4 of [22].
 By a similar method to that used in the proof of (i) in Theorem 3.2, one gets that (ii) and (iv) hold with the help of Theorem 4.3 of [22]. This completes the proof.\medskip

\noindent{\bf 3.3.  Inequalities among eigenvalues for  coupled BCs and those for some certain separated ones }\medskip

In this subsection,   we shall first establish inequalities among eigenvalues for a coupled BC  and those for some certain separated ones  applying  Theorems 3.2--3.5.  Then, for a fixed $K\in SL(2,\mathbb{R})$ and $\gamma\in(-\pi,0)\cup(0,\pi)$, we shall compare eigenvalues for $[K|-I]$,  those for $[e^{i\gamma}K|-I]$, and those for $[-K|-I]$. Combining the above two parts, we shall  establish inequalities among eigenvalues for three coupled BCs  and those for some certain separated ones, which generalize the main result of [18].\medskip

 Firstly, we shall establish inequalities among eigenvalues for a coupled BC  and those for some certain separated ones in the next two theorems.  Set $\lambda_n(e^{i\gamma}K):=\lambda_n([e^{i\gamma}K|-I])$ for briefness.\medskip

\noindent{\bf Theorem 3.6.} {\it Fix a difference equation $\pmb\omega=(1/f,q,w)$.  Let $\mathbf{A}=[e^{i\gamma}K|-I]$,  where $K\in SL(2,\mathbb{R})$ and $\gamma\in(-\pi,\pi]$.
 Then $(\pmb\omega,\mathbf{A})$ has exactly $N$ eigenvalues  if $k_{11}-f_0k_{12}\neq0$, and  exactly $N-1$ eigenvalues if $k_{11}-f_0k_{12}=0$;
 $(\pmb\omega,\mathbf{T}_K)$ has exactly $N$ eigenvalues if $k_{11}\neq0$, and  exactly $N-1$ eigenvalues if $k_{11}=0$;
  $(\pmb\omega,\mathbf{U}_K)$ has exactly $N-1$ eigenvalues  if $k_{11}-f_0k_{12}\neq0$, and exactly $N-2$ eigenvalues if $k_{11}-f_0k_{12}=0$, where \vspace{-0.2cm}$$\mathbf{T}_{K}:=\left [\begin{array} {cccc}0&1&0&0\\
0&0&-k_{21}&k_{11}\end{array}  \right ]\; and\;\;\mathbf{U}_K:=\left [\begin{array} {cccc}k_{11}&k_{12}&0&0\\
0&0&1&0\end{array}  \right ].\vspace{-0.2cm}$$
Furthermore, we have that
 \begin{itemize}\vspace{-0.2cm}
\item[{\rm (i)}] the eigenvalues of $(\pmb\omega,\mathbf{A})$ and  $(\pmb\omega,\mathbf{T}_K)$ satisfy the following inequalities:
\vspace{-0.2cm}$$\begin{array} {cccc} \lambda_0(\mathbf{T}_K)\leq \lambda_0(e^{i\gamma}K)\leq\lambda_1(\mathbf{T}_K)\leq\lambda_1(e^{i\gamma}K)\\[1.0ex]\leq\cdots\leq
\lambda_{N-1}(\mathbf{T}_K)\leq\lambda_{N-1}(e^{i\gamma}K) \end{array}\eqno(3.10)\vspace{-0.2cm}$$
 in the case that $(k_{11}-f_0k_{12})k_{11}f_0>0$;
\vspace{-0.2cm}$$\begin{array} {cccc}  \lambda_0(e^{i\gamma}K)\leq\lambda_0(\mathbf{T}_K)\leq\lambda_1(e^{i\gamma}K)\leq\lambda_1(\mathbf{T}_K)\\[1.0ex]\leq\cdots\leq
\lambda_{N-1}(e^{i\gamma}K)\leq\lambda_{N-1}(\mathbf{T}_K)\end{array}\eqno(3.11)\vspace{-0.2cm}$$
 in the case that $(k_{11}-f_0k_{12})k_{11}f_0<0$;
\vspace{-0.2cm}$$ \begin{array} {cccc} \lambda_0(\mathbf{T}_K)\leq  \lambda_0(e^{i\gamma}K)\leq\lambda_1(\mathbf{T}_K)\leq\lambda_1(e^{i\gamma}K)\\[1.0ex]\leq\cdots\leq
\lambda_{N-2}(\mathbf{T}_K)\leq\lambda_{N-2}(e^{i\gamma}K)\leq\lambda_{N-1}(\mathbf{T}_K)\end{array}\eqno(3.12)\vspace{-0.2cm}$$
 in the case that $k_{11}-f_0k_{12}=0$;
\vspace{-0.2cm}$$  \begin{array} {cccc}\lambda_0(e^{i\gamma}K)\leq \lambda_0(\mathbf{T}_K) \leq\lambda_1(e^{i\gamma}K)\leq\lambda_1(\mathbf{T}_K)\\[1.0ex]\leq\cdots\leq\lambda_{N-2}(e^{i\gamma}K)\leq
\lambda_{N-2}(\mathbf{T}_K)\leq\lambda_{N-1}(e^{i\gamma}K)\end{array}\eqno(3.13)
\vspace{-0.2cm}$$
in the case that $k_{11}=0$;
\item[{\rm (ii)}] the eigenvalues of $(\pmb\omega,\mathbf{A})$ and $(\pmb\omega,\mathbf{U}_K)$ satisfy the following inequalities:
\vspace{-0.2cm}$$  \begin{array} {cccc}\lambda_0(e^{i\gamma}K)\leq \lambda_0(\mathbf{U}_K) \leq\lambda_1(e^{i\gamma}K)\leq\lambda_1(\mathbf{U}_K)\\[1.0ex]\leq\cdots\leq\lambda_{N-2}(e^{i\gamma}K)\leq
\lambda_{N-2}(\mathbf{U}_K)\leq\lambda_{N-1}(e^{i\gamma}K)\end{array}\eqno(3.14)
\vspace{-0.2cm}$$
 in the case that  $k_{11}-f_0k_{12}\neq0$;
\vspace{-0.2cm} $$ \begin{array} {cccc}\lambda_0(e^{i\gamma}K)\leq  \lambda_0(\mathbf{U}_K) \leq\lambda_1(e^{i\gamma}K)\leq \lambda_1(\mathbf{U}_K)\\[1.0ex]\leq\cdots\leq\lambda_{N-3}(e^{i\gamma}K)\leq
 \lambda_{N-3}(\mathbf{U}_K)\leq\lambda_{N-2}(e^{i\gamma}K)\end{array}\eqno(3.15)
\vspace{-0.2cm}$$
 in the case that  $k_{11}-f_0k_{12}=0$.
 \end{itemize} }\medskip

\noindent{\bf Proof.} The number of eigenvalues of $(\pmb\omega,\mathbf{A})$, $(\pmb\omega,\mathbf{T}_K)$, and $(\pmb\omega,\mathbf{U}_K)$ in each case can be obtained by Lemma 2.4 and direct computations.
Let $k_{11}\neq 0$. Since $\det K =1$,
\vspace{-0.2cm}$$\begin{array} {cccc}\mathbf{A}=[e^{i\gamma}K|-I ] =\left[\begin{array} {cccc}
1&k_{12}/k_{11} & -e^{-i\gamma}/k_{11}&0\\
-e^{i\gamma}k_{21}&-e^{i\gamma}k_{22}&0&1\end{array}\right]
=\left[\begin{array} {cccc}
1&a_{12} &\bar{z}&0\\
0&z&b_{21}&1\end{array}\right]\in \mathcal{O}^{\mathbb{C}}_{1,4},\end{array}\vspace{-0.1cm}$$
where $a_{12}:=k_{12}/k_{11}$, $b_{21}:=-k_{21}/k_{11}$, and $z:=-e^{i\gamma}/k_{11}$.
Then by (i) of Lemma 2.7, $\mathbf{A}\in\mathcal{C}_{1,4,z, b_{21}}\cap\mathcal{C}_{1,4,z, a_{12}}$, and the corresponding LBCs  satisfy that
$$\vspace{-0.2cm}\begin{array} {cccc}\mathbf{S}_1=\left [\begin{array} {cccc}0&1&0&0\\
0&0&b_{21}&1\end{array}  \right ]=\left [\begin{array} {cccc}0&1&0&0\\
0&0&-k_{21}&k_{11}\end{array}  \right ]=\mathbf{T}_K,\\[2.0ex]
\mathbf{S}_2=\left [\begin{array} {cccc}1&a_{12}&0&0\\
0&0&1&0\end{array}  \right ]=\left [\begin{array} {cccc}k_{11}&k_{12}&0&0\\
0&0&1&0\end{array}  \right ]=\mathbf{U}_K.\end{array}\vspace{-0.05cm}$$
Note that $(k_{11}-f_0k_{12})k_{11}f_0>0$, $(k_{11}-f_0k_{12})k_{11}f_0<0$, $k_{11}-f_0k_{12}=0$, and $k_{11}-f_0k_{12}\neq0$  are equivalent to $a_{12}<1/f_0$, $a_{12}>1/f_0$,  $a_{12}=1/f_0$, and $a_{12}\neq1/f_0$, respectively.
Therefore, by Theorem 3.2, one gets that  $(k_{11}-f_0k_{12})k_{11}f_0>0$ implies (3.10); $(k_{11}-f_0k_{12})k_{11}f_0<0$ implies (3.11); $k_{11}-f_0k_{12}=0$ implies (3.12) and (3.15); $k_{11}-f_0k_{12}\neq0$ implies (3.14).

 Let $k_{11}=0$. Now we show that (3.13)--(3.14) hold in this case. Since  $k_{11}=0$, $-k_{12}k_{21}=1$. Denote
$$\vspace{-0.2cm} K_{\epsilon}:=\left (\begin{array} {cccc}\epsilon&k_{12}\\
(-1+\epsilon k_{22})/k_{12}&k_{22}\end{array}  \right )\in SL(2,\mathbb{R}),\;\;\; \;\epsilon\in\mathbb{R}.$$\vspace{-0.2cm}
Then $\lim\limits_{\epsilon\rightarrow0}K_{\epsilon}=K$.
By the definition of $\mathbf{T}_K$ and $\mathbf{U}_K$, we see that
$$\vspace{-0.2cm}\mathbf{T}_{K_{\epsilon}}=\left [\begin{array} {cccc}0&1&0&0\\
0&0&1-\epsilon k_{22}&\epsilon k_{12}\end{array}  \right ]\;{\rm and}\;\;\mathbf{U}_{K_{\epsilon}}=\left [\begin{array} {cccc}\epsilon&k_{12}&0&0\\
0&0&1&0\end{array}  \right ].$$\vspace{-0.2cm}
Then
\vspace{-0.2cm}$$[e^{i\gamma}K_{\epsilon}|-I]\rightarrow[e^{i\gamma}K|-I],\;\;\mathbf{T}_{K_{\epsilon}}\rightarrow\mathbf{T}_{K},\;\;
\mathbf{U}_{K_{\epsilon}}\rightarrow\mathbf{U}_{K},
\;{\rm as} \;\;\epsilon\rightarrow0.\vspace{0.05cm}$$

Since $k_{12}\neq0$,  one can choose a sufficiently small $\epsilon_1>0$ such that $\epsilon-f_0k_{12}\neq0$, where $0\leq \epsilon\leq \epsilon_1$. Thus, by Lemma 2.4, there are exactly $N$ eigenvalues for each $[e^{i\gamma}K_{\epsilon}|-I]$, $0\leq\epsilon\leq\epsilon_1$, and by Lemma 2.5, $\lambda_n(e^{i\gamma}K_{\epsilon})$ is continuous in $\epsilon\in[0,\epsilon_1]$,
which implies that
\vspace{-0.2cm}$$\lambda_n(e^{i\gamma}K_{\epsilon})\rightarrow\lambda_n(e^{i\gamma}K),\;{\rm as} \;\;\epsilon\rightarrow0^+,\;0\leq n\leq N-1.\eqno(3.16)\vspace{-0.2cm}$$

Suppose that $k_{12}<0$. By Lemma 2.4, $(\pmb\omega,\mathbf{U}_{K_\epsilon})$ has exactly $N-1$ eigenvalues for each $\epsilon\in[0,\epsilon_1]$. Thus by Lemma 2.5, $\lambda_n(\mathbf{U}_{K_\epsilon})$ is continuous in $\epsilon\in[0,\epsilon_1]$,
which implies that
\vspace{-0.2cm}$$\lambda_n(\mathbf{U}_{K_\epsilon})\rightarrow\lambda_n(\mathbf{U}_{K}),\;{\rm as} \;\;\epsilon\rightarrow0^+,\; 0\leq n\leq N-2.\eqno(3.17)\vspace{-0.2cm}$$
If $f_0>0$, then $\mathbf{T}_K\in\mathcal{B}_{2,3r} $ and $\mathbf{T}_{K_\epsilon}\in\mathcal{B}_{2,3r}^+ $, where $\epsilon\in(0,\epsilon_1]$ and
\vspace{-0.2cm}$$\begin{array}{cccc}
\mathcal{B}_{2,3r}:=\left\{ \mathbf{A}\in\mathcal{O}_{2,3}^{\mathbb{C}}:  (a_{11}+ f_0)b_{22}=|z|^2, a_{11}+  f_0\geq 0, b_{22}\geq 0\right\}\backslash\left\{\mathbf{C}\right\},\\[1.0ex]
\mathcal{B}_{2,3r}^+:=\left\{ \mathbf{A}\in\mathcal{O}_{2,3}^{\mathbb{C}}: a_{11}\geq- f_0,b_{22}\geq 0, (a_{11}+ f_0)b_{22}>|z|^2\right\}.\end{array}\vspace{-0.1cm}$$
Note that $\epsilon_1$ can be chosen such that $1-\epsilon k_{22}>0$ for any $0<\epsilon\leq\epsilon_1$.
By Theorem 4.4 in [22],
\vspace{-0.2cm}$$\lambda_0(\mathbf{T}_{K_\epsilon})\rightarrow-\infty,\; \lambda_{n}(\mathbf{T}_{K_\epsilon})\rightarrow\lambda_{n-1}(\mathbf{T}_{K}),\; as\; \epsilon\rightarrow0^+,\;1\leq n\leq N-1.\eqno(3.18)\vspace{-0.2cm}$$
If $f_0<0$, then $\mathbf{T}_K\in\mathcal{B}_{2,3l}$ and $\mathbf{T}_{K_\epsilon}\in\mathcal{B}_{2,3}^-$, where $\epsilon\in(0,\epsilon_1]$ and
\vspace{-0.2cm}$$ \begin{array}{cccc}\mathcal{B}_{2,3l}:=\left\{ \mathbf{A}\in\mathcal{O}_{2,3}^{\mathbb{C}}:  (a_{11}+ f_0)b_{22}=|z|^2, a_{11}+  f_0\leq 0, b_{22}\leq 0\right\}
\backslash\left\{\mathbf{C}\right\},\\[1.0ex]
\mathcal{B}_{2,3}^-:=\left\{ \mathbf{A}\in\mathcal{O}_{2,3}^{\mathbb{C}}:  (a_{11}+ f_0)b_{22}<|z|^2\right\}.\end{array}\vspace{-0.01cm}$$
 By Theorem 4.4 in [22], (3.18) holds.

Since $(\epsilon-f_0k_{12})f_0\epsilon>0$, where $0<\epsilon\leq \epsilon_1$, by (3.10) and (3.14) for $[e^{i\gamma}K_{\epsilon}|-I]$,
\vspace{-0.2cm}$$\begin{array} {cccc}  \lambda_0(\mathbf{T}_{K_\epsilon})\leq\lambda_0(e^{i\gamma}K_\epsilon)\leq \{\lambda_1(\mathbf{T}_{K_\epsilon}),\lambda_0(\mathbf{U}_{K_\epsilon})\} \leq\lambda_1(e^{i\gamma}K_\epsilon)\leq\{\lambda_2(\mathbf{T}_{K_\epsilon}),\\[1.0ex]\lambda_1(\mathbf{U}_{K_\epsilon})\}\leq\cdots\leq
\lambda_{N-2}(e^{i\gamma}K_\epsilon)\leq\{\lambda_{N-1}(\mathbf{T}_{K_\epsilon}),\lambda_{N-2}(\mathbf{U}_{K_\epsilon})\}\leq\lambda_{N-1}(e^{i\gamma}K_\epsilon).
\end{array}\eqno(3.19)\vspace{-0.2cm}$$
Let $\epsilon\rightarrow0^+$ in (3.19), it follows from (3.16)--(3.18) that (3.13)--(3.14) hold  for $[e^{i\gamma}K|-I]$.

Suppose that $k_{12}>0$. With a similar method to that used  in the case that
$k_{12}<0$, one can show that   (3.13)--(3.14) hold for $[e^{i\gamma}K|-I]$. The proof is   complete.\medskip

\noindent{\bf Theorem 3.7.} {\it  Fix a difference equation $\pmb\omega=(1/f,q,w)$.  Let $\mathbf{A}=[e^{i\gamma}K|-I]$, where $K\in SL(2,\mathbb{R})$ and $\gamma\in(-\pi,\pi]$.
 Then $(\pmb\omega, \mathbf{S}_K)$ has exactly $N$ eigenvalues if $k_{12}\neq0$, and  exactly $N-1$ eigenvalues if $k_{12}=0$; $(\pmb\omega, \mathbf{V}_K)$ has exactly $N$ eigenvalues if $f_0k_{22}-k_{21}\neq0$, and exactly $N-1$ eigenvalues  if $f_0k_{22}-k_{21}=0$, where \vspace{-0.2cm}$$\mathbf{S}_{K}:=\left [\begin{array} {cccc}1&0&0&0\\
0&0&-k_{22}&k_{12}\end{array}  \right ]\;and \;\; \mathbf{V}_K:=\left [\begin{array} {cccc}k_{21}&k_{22}&0&0\\
0&0&0&1\end{array}  \right ].\vspace{-0.2cm}$$  Furthermore, we have that
\begin{itemize}\vspace{-0.2cm}
\item[{\rm (i)}] the eigenvalues of $(\pmb\omega,\mathbf{A})$ and $(\pmb\omega,\mathbf{S}_K)$ satisfy the following inequalities:
\vspace{-0.2cm} $$\begin{array}  {cccc}
\lambda_0(e^{i\gamma}K)\leq\lambda_{0}(\mathbf{S}_K)\leq\lambda_{1}(e^{i\gamma}K)\leq\lambda_{1}(\mathbf{S}_K)\\[1.0ex]\leq\cdots\leq
\lambda_{N-2}(e^{i\gamma}K)\leq\lambda_{N-2}(\mathbf{S}_K)\leq \lambda_{N-1}(e^{i\gamma}K)\end{array} \eqno(3.20)\vspace{-0.2cm} $$
in the case that  $k_{12}=0$;
\vspace{-0.2cm} $$\begin{array} {cccc}
\lambda_0(\mathbf{S}_K)\leq\lambda_0(e^{i\gamma}K)\leq\lambda_1(\mathbf{S}_K)\leq\lambda_1(e^{i\gamma}K)\\[1.0ex]\leq\cdots\leq\lambda_{N-1}(\mathbf{S}_K)\leq \lambda_{N-1}(e^{i\gamma}K)
\end{array}\eqno(3.21)\vspace{-0.2cm} $$
in the case that $(k_{11}-f_0k_{12})k_{12}>0$;
\vspace{-0.2cm} $$\begin{array} {cccc}
\lambda_0(e^{i\gamma}K)\leq\lambda_0(\mathbf{S}_K)\leq\lambda_1(e^{i\gamma}K)\leq\lambda_1(\mathbf{S}_K)\\[1.0ex]\leq\cdots
\leq\lambda_{N-1}(e^{i\gamma}K)\leq \lambda_{N-1}(\mathbf{S}_K)\end{array}\eqno(3.22)\vspace{-0.2cm} $$
in the case that $(k_{11}-f_0k_{12})k_{12}<0$;
\vspace{-0.2cm} $$\begin{array} {cccc}
\lambda_0(\mathbf{S}_K)\leq\lambda_0(e^{i\gamma}K)\leq\lambda_1(\mathbf{S}_K)\leq\lambda_1(e^{i\gamma}K)\\[1.0ex]\leq\cdots
\leq\lambda_{N-2}(\mathbf{S}_K)\leq\lambda_{N-2}(e^{i\gamma}K)\leq\lambda_{N-1}(\mathbf{S}_K)
\end{array}\eqno(3.23)\vspace{-0.2cm} $$
in the case that $k_{11}-f_0k_{12}=0$;
\item[{\rm (ii)}] the eigenvalues of $(\pmb\omega,\mathbf{A})$ and $(\pmb\omega,\mathbf{V}_K)$ satisfy the following inequalities:
\vspace{-0.2cm} $$\begin{array} {cccc}
\lambda_0(e^{i\gamma}K)\leq\lambda_0(\mathbf{V}_K)\leq\lambda_1(e^{i\gamma}K)\leq\lambda_1(\mathbf{V}_K)\\[1.0ex]\leq\cdots\leq\lambda_{N-2}(e^{i\gamma}K)\leq \lambda_{N-2}(\mathbf{V}_K)\leq\lambda_{N-1}(e^{i\gamma}K)
\end{array}\eqno(3.24)\vspace{-0.2cm} $$
in the case that  $f_0k_{22}-k_{21}=0$;
\vspace{-0.2cm} $$\begin{array} {cccc}
\lambda_0(\mathbf{V}_K)\leq\lambda_0(e^{i\gamma}K)\leq\lambda_1(\mathbf{V}_K)\leq\lambda_1(e^{i\gamma}K)\\[1.0ex]\leq\cdots\leq \lambda_{N-1}(\mathbf{V}_K)\leq\lambda_{N-1}(e^{i\gamma}K)
\end{array}\eqno(3.25)\vspace{-0.2cm} $$
in the case that $(k_{11}-f_0k_{12})(f_0k_{22}-k_{21})>0$;
\vspace{-0.2cm} $$\begin{array} {cccc}
\lambda_0(e^{i\gamma}K)\leq\lambda_0(\mathbf{V}_K)\leq\lambda_1(e^{i\gamma}K)\leq\lambda_1(\mathbf{V}_K)\\[1.0ex]\leq\cdots\leq \lambda_{N-1}(e^{i\gamma}K)\leq\lambda_{N-1}(\mathbf{V}_K)
\end{array}\eqno(3.26)\vspace{-0.2cm} $$
in the case that $(k_{11}-f_0k_{12})(f_0k_{22}-k_{21})<0$;
\vspace{-0.2cm} $$\begin{array} {cccc}
\lambda_0(\mathbf{V}_K)\leq\lambda_0(e^{i\gamma}K)\leq\lambda_1(\mathbf{V}_K)\leq\lambda_1(e^{i\gamma}K)\\[1.0ex]\leq\cdots\leq\lambda_{N-2}(\mathbf{V}_K)\leq \lambda_{N-2}(e^{i\gamma}K)\leq\lambda_{N-1}(\mathbf{V}_K)
\end{array}\eqno(3.27)\vspace{-0.2cm} $$
in the case that $k_{11}-f_0k_{12}=0$.
\end{itemize}\vspace{-0.2cm}}

\noindent{\bf Proof.} The number of eigenvalues of $(\pmb\omega,\mathbf{A})$, $(\pmb\omega,\mathbf{S}_K)$, and $(\pmb\omega,\mathbf{V}_K)$ in each case can be obtained by Lemma 2.4 and direct computations.
Let $k_{22}\neq0$.  Since $\det K =1$,
$$\vspace{-0.2cm}\mathbf{A}=[e^{i\gamma}K|-I ] =\left[\begin{array} {llll}
e^{i\gamma}k_{11}&e^{i\gamma}k_{12}&-1&0 \\
-k_{21}/k_{22}&-1&0&e^{-i\gamma}/k_{22}\end{array}\right]=
\left[\begin{array} {llll}
a_{11}&-1&0&\bar{z} \\
z&0&-1&b_{22}\end{array}\right]\in \mathcal{O}^{\mathbb{C}}_{2,3},\vspace{0.1cm}$$
where $a_{11}:=-k_{21}/k_{22}$, $z:=e^{i\gamma}/k_{22}$, $b_{22}:=k_{12}/k_{22}$.
Then by (iii) of Lemma 2.7, $\mathbf{A}\in\mathcal{C}_{2,3,z, b_{22}}\cap\mathcal{C}_{2,3,z, a_{11}}$, and the corresponding LBCs  satisfy that
$$\vspace{-0.2cm}\begin{array} {cccc}\mathbf{S}_5=\left [\begin{array} {cccc}1&0&0&0\\
0&0&-1&b_{22}\end{array}  \right ]=\left [\begin{array} {cccc}1&0&0&0\\
0&0&-k_{22}&k_{12}\end{array}  \right ]=\mathbf{S}_K,\\[2.0ex]
\mathbf{S}_6=\left [\begin{array} {cccc}a_{11}&-1&0&0\\
0&0&0&1\end{array}  \right ]=\left [\begin{array} {cccc}k_{21}&k_{22}&0&0\\
0&0&0&1\end{array}  \right ]=\mathbf{V}_K.\end{array}\vspace{0.05cm}$$
Note that $b_{22}=0$ is equivalent to $k_{12}=0$; in the case that $b_{22}\neq0$, one gets that   $a_{11}+f_0<|z|^2/b_{22}$,  $a_{11}+f_0>|z|^2/b_{22}$, and $a_{11}+f_0=|z|^2/b_{22}$ are equivalent to $(k_{11}-f_0k_{12})k_{12}>0$,
 $(k_{11}-f_0k_{12})k_{12}<0$,
and  $k_{11}-f_0k_{12}=0$, respectively;
 $a_{11}+f_0=0$ is equivalent to $f_0k_{22}-k_{21}=0$;
in the case that $a_{11}+f_0\neq0$, one gets that  $b_{22}<|z|^2/(a_{11}+f_0)$, $b_{22}>|z|^2/(a_{11}+f_0)$, and $b_{22}=|z|^2/(a_{11}+f_0)$
 are equivalent to $(k_{11}-f_0k_{12})(f_0k_{22}-k_{21})>0$,   $(k_{11}-f_0k_{12})(f_0k_{22}-k_{21})<0$, and $k_{11}-f_0k_{12}=0$, respectively.
Therefore, by Theorem 3.4, one gets that   $k_{12}=0$ implies (3.20); $(k_{11}-f_0k_{12})k_{12}>0$ implies (3.21); $(k_{11}-f_0k_{12})k_{12}<0$ implies (3.22); $k_{11}-f_0k_{12}=0$ implies (3.23) and (3.27); $f_0k_{22}-k_{21}=0$ implies (3.24); $(k_{11}-f_0k_{12})(f_0k_{22}-k_{21})>0$ implies (3.25);
$(k_{11}-f_0k_{12})(f_0k_{22}-k_{21})<0$ implies (3.26).

Let $k_{22}=0$.  Now we show that (3.21)--(3.23) and (3.25)--(3.27) hold in this case. Since $k_{22}=0$, $-k_{12}k_{21}=1$.
Denote
$$\vspace{-0.2cm} K_{\epsilon}=\left (\begin{array} {cccc}k_{11}&k_{12}\\
(-1+\epsilon k_{11})/k_{12}&\epsilon \end{array}  \right )\in SL(2,\mathbb{R}), \;\;\epsilon\in\mathbb{R}.\eqno(3.28)\vspace{-0.1cm}$$
Then $\lim\limits_{\epsilon\rightarrow0}K_\epsilon=K$. By the definition of $\mathbf{S}_K$ and $\mathbf{V}_K$, one has that
\vspace{-0.2cm}$$ \mathbf{S}_{K_{\epsilon}}=\left [\begin{array} {cccc}1&0&0&0\\
0&0&-\epsilon&k_{12}\end{array}  \right ]\;\;{\rm and} \;\;\begin{array} {cccc}\mathbf{V}_{K_\epsilon}=\left [\begin{array} {cccc}-1+\epsilon k_{11}& \epsilon k_{12}&0&0\\
0&0&0&1\end{array}  \right ].\end{array}\vspace{-0.2cm}$$
Then
$$\vspace{-0.1cm}[e^{i\gamma}K_{\epsilon}|-I]\rightarrow[e^{i\gamma}K|-I],\;\;\mathbf{S}_{K_{\epsilon}}\rightarrow\mathbf{S}_{K},\;\;
\mathbf{V}_{K_{\epsilon}}\rightarrow\mathbf{V}_{K},
\;{\rm as} \;\;\epsilon\rightarrow0.\vspace{0.05cm}$$

In the case that  $(k_{11}-f_0k_{12})k_{12}>0,$  by (3.21) for $[e^{i\gamma}K_\epsilon|-I]$, $\epsilon>0$, one gets that
\vspace{-0.2cm} $$\begin{array} {cccc}
\lambda_0(\mathbf{S}_{K_{\epsilon}})\leq\lambda_0(e^{i\gamma}K_\epsilon)\leq\lambda_1(\mathbf{S}_{K_{\epsilon}})\leq\lambda_1(e^{i\gamma}K_\epsilon)\\[1.0ex]\leq\cdots
\leq\lambda_{N-1}(\mathbf{S}_{K_{\epsilon}})\leq\lambda_{N-1}(e^{i\gamma}K_\epsilon). \end{array}\eqno(3.29)\vspace{-0.2cm} $$
Since $k_{11}-f_0k_{12}\neq0$ and $k_{12}\neq0$,  there are exactly $N$ eigenvalues for each $[e^{i\gamma}K_\epsilon|-I]$ and for each $\mathbf{S}_{K_\epsilon}$, where $0\leq\epsilon\leq1$, by Lemma 2.4. It follows from  Lemma 2.5 that $\lambda_n(e^{i\gamma}K_\epsilon)$ and $\lambda_n(\mathbf{S}_{K_\epsilon})$ are continuous in $\epsilon\in[0,1]$, which implies that,
\vspace{-0.2cm} $$\lambda_n(e^{i\gamma}K_\epsilon)\rightarrow\lambda_n(e^{i\gamma}K),\;\lambda_n(\mathbf{S}_{K_\epsilon})\rightarrow\lambda_n(\mathbf{S}_{K}), \;\;{\rm as} \;\; \epsilon\rightarrow0^+,\;0\leq n \leq N-1.\eqno(3.30)\vspace{-0.2cm}$$
 Let $\epsilon\rightarrow0^+$ in (3.29), it follows from (3.30) that (3.21) holds for $[e^{i\gamma}K|-I]$.

 With  similar arguments to the proof of (3.21) for $[e^{i\gamma}K|-I]$, one can show that (3.22)--(3.23) hold for $[e^{i\gamma}K|-I]$.

Next, we show that (3.25)--(3.27) hold for $[e^{i\gamma}K|-I]$. In the case that  $-(k_{11}-f_0k_{12})k_{21}>0$,  one can choose an $\epsilon_1>0$ sufficiently small that
 $(k_{11}-f_0k_{12})(f_0\epsilon-(-1+\epsilon k_{11})/k_{12})=-(k_{11}-f_0k_{12})(k_{21}+\epsilon(k_{11}-f_0k_{12})/k_{12})>0$ and $1-\epsilon(k_{11}-f_0k_{12})>0$, $0\leq\epsilon\leq\epsilon_1$. Then by $(3.25)$ for $[e^{i\gamma}K_{\epsilon}|-I]$, where $0< \epsilon\leq\epsilon_1$, one has that
 \vspace{-0.2cm} $$\begin{array} {cccc}
\lambda_0(\mathbf{V}_{K_\epsilon})\leq\lambda_0(e^{i\gamma}K_\epsilon)\leq\lambda_1(\mathbf{V}_{K_\epsilon})\leq\lambda_1(e^{i\gamma}K_\epsilon)\\[1.0ex]\leq\cdots\leq \lambda_{N-1}(\mathbf{V}_{K_\epsilon})\leq\lambda_{N-1}(e^{i\gamma}K_\epsilon).
\end{array}\eqno(3.31)\vspace{-0.2cm} $$
 Since $k_{11}-f_0k_{12}\neq0$ and  $1- k_{11}\epsilon+f_0 k_{12}\epsilon=1-\epsilon(k_{11}-f_0k_{12})>0$, by Lemma 2.4 there are exactly $N$ eigenvalues for each $[e^{i\gamma}K_\epsilon|-I]$ and each $\mathbf{V}_{K_\epsilon}$, where $0\leq\epsilon\leq\epsilon_1$. By Lemma 2.5, $\lambda_n(e^{i\gamma}K_\epsilon)$ and $\lambda_n(\mathbf{V}_{K_\epsilon})$ are continuous in $\epsilon\in[0,\epsilon_1]$, which implies that
  \vspace{-0.2cm}$$\lambda_n(e^{i\gamma}K_\epsilon)\rightarrow\lambda_n(e^{i\gamma}K),\;\;\lambda_n(\mathbf{V}_{K_\epsilon})\rightarrow\lambda_n(\mathbf{V}_{K}),\;\;{\rm as}\;\;\epsilon\rightarrow0^+,\;0\leq n\leq N-1.\eqno(3.32)
  \vspace{-0.2cm}$$
  Let $\epsilon\rightarrow0^+$ in $(3.31)$, it follows from $(3.32)$ that $(3.25)$ holds for $[e^{i\gamma}K|-I]$.

 With a similar argument to the proof of (3.25) for $[e^{i\gamma}K|-I]$, one can show that (3.26)--(3.27) hold for $[e^{i\gamma}K|-I]$. This completes the proof.
 \medskip

\noindent{\bf Remark 3.2.} (ii) of  Theorem 3.6 and (i) of Theorem 3.7 can  also be obtained by dividing the  discussion into two cases: $k_{12}\neq0$ and $k_{12}=0$,  applying Theorem 3.3, and using  a similar method to that used in  the  proof of them; while (i) of Theorem 3.6 and (ii) of Theorem 3.7 can  also be obtained by dividing the  discussion into two cases: $k_{21}\neq0$ and $k_{21}=0$, applying Theorem 3.5, and using   a similar method to that used in  the  proof of them.\medskip

The following result, which is a direct consequence of Theorems 3.6--3.7, gives comparison of eigenvalues for $[e^{i\gamma}K|-I]$ with those for $\mathbf{S}_K$, those for $\mathbf{U}_K$, those for $\mathbf{T}_K$, and those for $\mathbf{V}_K$ under the assumption that $k_{11}-f_0k_{12}=0$.\medskip

\noindent{\bf Corollary 3.2.} {\it Fix a difference equation $\pmb\omega=(1/f,q,w)$. Let $\mathbf{A}=[e^{i\gamma}K|-I]\in \mathcal{B}^{\mathbb{C}}$, where $K\in SL(2,\mathbb{R})$ and $\gamma\in(-\pi,\pi]$. If $k_{11}-f_0k_{12}=0$, then
\vspace{-0.2cm} $$\begin{array} {cccc}
\{\lambda_0(\mathbf{S}_K),\lambda_0(\mathbf{T}_K),\lambda_0(\mathbf{V}_K)\}\leq\lambda_0(e^{i\gamma}K)\leq\{\lambda_1(\mathbf{S}_K),\lambda_1(\mathbf{T}_K),
\lambda_1(\mathbf{V}_K),\lambda_0(\mathbf{U}_K)\}\leq\\[1.0ex]\lambda_1(e^{i\gamma}K)\leq
\{\lambda_2(\mathbf{S}_K),\lambda_2(\mathbf{T}_K),
\lambda_2(\mathbf{V}_K),\lambda_1(\mathbf{U}_K)\}\leq\cdots\leq\lambda_{N-3}(e^{i\gamma}K)\leq\{\lambda_{N-2}(\mathbf{S}_K),\\[1.0ex]\lambda_{N-2}(\mathbf{T}_K),
\lambda_{N-2}(\mathbf{V}_K),\lambda_{N-3}(\mathbf{U}_K)\}
\leq\lambda_{N-2}(e^{i\gamma}K)\leq\{\lambda_{N-1}(\mathbf{S}_K),\lambda_{N-1}(\mathbf{T}_K),
\lambda_{N-1}(\mathbf{V}_K)\}.
\end{array} $$\vspace{-0.2cm}}

Note that a coupled BC $[e^{i\gamma}K|-I]$ can be written as  $[e^{i\gamma/2}K|-e^{-i\gamma/2}I]$.  Then by Lemma 2.3, a simple calculation yields that
\vspace{-0.2cm}$$\Gamma(\lambda)=2\cos\gamma-k_{22}\phi_{N}(\lambda)+k_{21}\psi_{N}(\lambda)+k_{12}f_N\Delta\phi_N(\lambda)-k_{11}f_N\Delta\psi_N(\lambda).
\eqno(3.33)
\vspace{-0.2cm}$$
Thus,  the eigenvalues for $[e^{i\gamma}K|-I]$ are the same as  those for $[e^{-i\gamma}K|-I]$ by (3.33).
Now, it's ready to establish inequalities among eigenvalues for the three coupled BCs: $[K|-I]$, $[e^{i\gamma}K|-I]$, and $[-K|-I]$, and those for  the corresponding separated ones.
\medskip

\noindent{\bf Theorem 3.8.} {\it Fix a difference equation $\pmb\omega=(1/f,q,w)$ satisfying that $\prod_{i=0}^{N-1}(1/f_i)>0.$  Let $\gamma\in(-\pi,0)\cup(0,\pi)$ and $K\in SL(2,\mathbb{R})$ satisfy that $k_{11}-f_0k_{12}\neq0$. Then the eigenvalues of SLPs $(\pmb\omega, [K|-I])$,
$(\pmb\omega, [e^{i\gamma}K|-I])$, $(\pmb\omega, [-K|-I])$, and $(\pmb\omega, \mathbf{S}_K)$ satisfy the following inequalities:
 \begin{itemize}\vspace{-0.2cm}
\item[{\rm (i)}] for $k_{11}-f_0k_{12}>0$ and $k_{12}>0$,
\vspace{-0.1cm} $$\begin{array} {cccc}\lambda_0(\mathbf{S}_K)\leq\lambda_0(K)<\lambda_0(e^{i\gamma}K)<\lambda_0(-K)\leq\lambda_1(\mathbf{S}_K)\leq
\lambda_1(-K)<\\[1.0ex]\lambda_1(e^{i\gamma}K)<\lambda_1(K)
 \leq\cdots\leq\lambda_{N-2}(\mathbf{S}_K)\leq\lambda_{N-2}(K)<\lambda_{N-2}(e^{i\gamma}K)<\\[1.0ex]\lambda_{N-2}(-K)
 \leq\lambda_{N-1}(\mathbf{S}_K)\leq\lambda_{N-1}(-K)
 <\lambda_{N-1}(e^{i\gamma}K)<\lambda_{N-1}(K)\end{array}\eqno(3.34)\vspace{-0.2cm}$$
 in the case that $N$ is even;
\vspace{-0.1cm} $$\begin{array} {cccc}\lambda_0(\mathbf{S}_K)\leq\lambda_0(K)<\lambda_0(e^{i\gamma}K)<\lambda_0(-K)\leq\lambda_1(\mathbf{S}_K)\leq
\lambda_1(-K)<\\[1.0ex]\lambda_1(e^{i\gamma}K)<\lambda_1(K)
 \leq\cdots\leq\lambda_{N-2}(\mathbf{S}_K)\leq\lambda_{N-2}(-K)<\lambda_{N-2}(e^{i\gamma}K)<\\[1.0ex]\lambda_{N-2}(K)
 \leq\lambda_{N-1}(\mathbf{S}_K)\leq\lambda_{N-1}(K)
 <\lambda_{N-1}(e^{i\gamma}K)<\lambda_{N-1}(-K)\end{array}\eqno(3.35)\vspace{-0.2cm}$$
  in the case that $N$ is odd;
  \item[{\rm (ii)}] for $k_{11}-f_0k_{
  12}>0$ and $k_{12}<0$,
  \vspace{-0.1cm} $$\begin{array} {cccc}\lambda_0(K)<\lambda_0(e^{i\gamma}K)<\lambda_0(-K)\leq\lambda_0(\mathbf{S}_K)\leq
\lambda_1(-K)<\lambda_1(e^{i\gamma}K)<\lambda_1(K)\\[1.0ex]\leq\lambda_1(\mathbf{S}_K)
 \leq\cdots\leq\lambda_{N-2}(K)<\lambda_{N-2}(e^{i\gamma}K)<\lambda_{N-2}(-K)\leq\\[1.0ex]\lambda_{N-2}(\mathbf{S}_K)
 \leq\lambda_{N-1}(-K)<\lambda_{N-1}(e^{i\gamma}K)<\lambda_{N-1}(K)
 \leq\lambda_{N-1}(\mathbf{S}_K)
 \end{array}\vspace{-0.2cm}$$
 in the case that $N$ is even;
 \vspace{-0.1cm} $$\begin{array} {cccc}\lambda_0(K)<\lambda_0(e^{i\gamma}K)<\lambda_0(-K)\leq\lambda_0(\mathbf{S}_K)\leq
\lambda_1(-K)<\lambda_1(e^{i\gamma}K)<\lambda_1(K)\\[1.0ex]\leq\lambda_1(\mathbf{S}_K)
 \leq\cdots\leq\lambda_{N-2}(-K)<\lambda_{N-2}(e^{i\gamma}K)<\lambda_{N-2}(K)\leq\\[1.0ex]\lambda_{N-2}(\mathbf{S}_K)
 \leq\lambda_{N-1}(K)<\lambda_{N-1}(e^{i\gamma}K)<\lambda_{N-1}(-K)
 \leq\lambda_{N-1}(\mathbf{S}_K)\end{array}\vspace{-0.2cm}$$
  in the case that $N$ is odd;
  \item[{\rm (iii)}] for $k_{11}>0$ and $k_{12}=0$,
  \vspace{-0.1cm} $$\begin{array} {cccc}\lambda_0(K)<\lambda_0(e^{i\gamma}K)<\lambda_0(-K)\leq\lambda_0(\mathbf{S}_K)\leq
\lambda_1(-K)<\lambda_1(e^{i\gamma}K)\\[1.0ex]<\lambda_1(K)\leq\lambda_1(\mathbf{S}_K)
 \leq\cdots\leq\lambda_{N-2}(K)<\lambda_{N-2}(e^{i\gamma}K)<\lambda_{N-2}(-K)\\[1.0ex]\leq\lambda_{N-2}(\mathbf{S}_K)
 \leq\lambda_{N-1}(-K)<\lambda_{N-1}(e^{i\gamma}K)<\lambda_{N-1}(K)\end{array}\vspace{-0.2cm}$$
 in the case that $N$ is even;
  \vspace{-0.1cm} $$\begin{array} {cccc}\lambda_0(K)<\lambda_0(e^{i\gamma}K)<\lambda_0(-K)\leq\lambda_0(\mathbf{S}_K)\leq
\lambda_1(-K)<\lambda_1(e^{i\gamma}K)\\[1.0ex]<\lambda_1(K)\leq\lambda_1(\mathbf{S}_K)
 \leq\cdots\leq\lambda_{N-2}(-K)<\lambda_{N-2}(e^{i\gamma}K)<\lambda_{N-2}(K)\\[1.0ex]\leq\lambda_{N-2}(\mathbf{S}_K)
 \leq\lambda_{N-1}(K)<\lambda_{N-1}(e^{i\gamma}K)<\lambda_{N-1}(-K)\end{array}\vspace{-0.2cm}$$
  in the case that $N$ is odd;
   \item[{\rm (iv)}] assertions in  {\rm (i)--(iii)} hold with $K$ replaced by $-K$.
  \end{itemize}}

\noindent{\bf Proof.} First, we show that (i) holds. We only show that (3.34) holds, since (3.35) can be shown similarly. By  (2.3)  and (3.33), one can easily verify that the leading term of $\Gamma(\lambda)$ as a polynomial of $\lambda$ is
  \vspace{-0.2cm}$$ (-1)^{N+1}1/f_0(w_N\prod_{i=1}^{N-1}(w_i/f_i))(k_{11}-f_0k_{12})\lambda^N.  \vspace{-0.2cm}$$
  Since $k_{11}-f_0k_{12}>0$ and $1/f_0(w_N\prod_{i=1}^{N-1}(w_i/f_i))>0$, one has that
    \vspace{-0.2cm}$$ \lim\limits_{\lambda\rightarrow-\infty}\Gamma(\lambda)=-\infty,\;\;\lim\limits_{\lambda\rightarrow+\infty}\Gamma(\lambda)=-\infty,
        \eqno(3.36)\vspace{-0.2cm}$$
in the case that $N$ is even.

Let $\gamma\in(-\pi,0)\cup(0,\pi)$. By Lemma 2.3,
$\lambda_n(e^{i\gamma}K),$ $0\leq n\leq N-1$, are
exactly the zeros of the polynomial
$\Gamma(\lambda)$ for $(\pmb\omega, [e^{i\gamma}K|-I])$.  It follows from Theorem 3.1 in [23] that $\lambda_n(e^{i\gamma}K)$ is a simple eigenvalue for each $0\leq n\leq N-1$. Thus by Rolle mean value Theorem, there are exactly $N-1$ real  zeros  for  $\Gamma'(\lambda)$ and they are denoted by $x_1,\cdots, x_{N-1}$. Then  $x_n\in(\lambda_{n-1}(e^{i\gamma}K),\lambda_n(e^{i\gamma}K))$, and $\Gamma(\lambda)$ is strictly increasing in $(-\infty,x_1)$ and strictly decreasing in $(x_1,x_2)$. Hence $\Gamma(x_1)>0$. From (3.33), (3.36), and the above discussion, it follows that $ \lambda_0(K)<\lambda_0(e^{i\gamma}K)<\lambda_0(-K)$ and
$\lambda_1(-K)<\lambda_1(e^{i\gamma}K)<\lambda_1(K)$. Similarly, one can show that
\vspace{-0.2cm}$$\begin{array} {llll} \lambda_j(K)<\lambda_j(e^{i\gamma}K)<\lambda_j(-K),\;\;j=0,2,4,\cdots, N-2,\\[1.0ex]
 \lambda_j(-K)<\lambda_j(e^{i\gamma}K)<\lambda_j(K),\;\;j=1,3,5,\cdots, N-1.\end{array}\eqno(3.37) $$\vspace{-0.2cm}

Since $(k_{11}-f_0k_{12})k_{12}>0$, it follows from $(3.21)$ that
\vspace{-0.2cm}$$\begin{array} {cccc}\lambda_0(\mathbf{S}_K)\leq \{\lambda_0(e^{i\gamma}K):\gamma\in(-\pi,\pi]\}\leq\lambda_1(\mathbf{S}_K)\leq\{\lambda_1(e^{i\gamma}K):\gamma\in(-\pi,\pi]\}\\[1.0ex]\leq\cdots\leq
\lambda_{N-1}(\mathbf{S}_K)\leq\{\lambda_{N-1}(e^{i\gamma}K):\gamma\in(-\pi,\pi]\}. \end{array}\eqno(3.38)\vspace{-0.2cm}$$
Therefore, (3.37)--(3.38) imply that (3.34) holds.

Assertions (ii)--(iii) can be shown similarly.

Now we show that (iv) holds. It follows from the definition of $\mathbf{S}_K$ in Theorem 3.7 that $\mathbf{S}_K=\mathbf{S}_{-K}$.
Since $k_{11}-f_0k_{12}\neq0$ and the entries of $K$ satisfy none of the conditions in (i), (ii) or (iii), there are exactly three cases:
 \begin{itemize}\vspace{-0.2cm}
\item[{\rm (i$'$)}] $k_{11}-f_0k_{12}<0$, $k_{12}<0$;\vspace{-0.2cm}
\item[{\rm (ii$'$)}] $k_{11}-f_0k_{12}<0$, $k_{12}>0$;\vspace{-0.2cm}
\item[{\rm (iii$'$)}] $k_{11}<0$, $k_{12}=0$.\vspace{-0.2cm}
 \end{itemize}\vspace{-0.01cm}
If the entries of $K$ satisfy (i$'$), (ii$'$), and (iii$'$), separately, then assertions in (i), (ii), and (iii) hold for $-K$, respectively.
This completes the proof.\medskip

The following result is a direct consequence of Theorems 3.6 and 3.8.\medskip

\noindent{\bf Corollary 3.3.} {\it  Fix a difference equation $\pmb\omega=(1/f,q,w)$ satisfying that  $\prod_{i=0}^{N-1}(1/f_i)>0$. Let  $K\in SL(2,\mathbb{R})$.
\begin{itemize}\vspace{-0.3cm}
\item[{\rm (i)}] If $k_{11}-f_0k_{12}>0$, then
 there are exactly $N$ eigenvalues for $[K|-I]$ and $[-K|-I]$, separately. Further,  whether $N$ is odd or even, $\lambda_0(K)$  is a simple eigenvalue;  if $\lambda_j(K)<\lambda_{j+1}(K)$ for some odd number $j$ $(1\leq j\leq N-2)$, then $\lambda_j(K)$ and $\lambda_{j+1}(K)$ are
simple eigenvalues. Similar results hold in the case that $\lambda_j(-K)<\lambda_{j+1}(-K)$ for some even  number $j$ $(0\leq j\leq N-2)$.
If $N$ is odd, then $\lambda_{N-1}(-K)$ is a simple eigenvalue; and if $N$ is
even, then $\lambda_{N-1}(K)$ is a simple eigenvalue.\vspace{-0.3cm}
\item[{\rm (ii)}]If $k_{11}-f_0k_{12}<0$,  similar results in {\rm (i)} can be obtained with $K$  replaced by $-K$.
\end{itemize}}\medskip

\noindent{\bf Remark 3.3.} Theorem 3.1 in [18] gives inequalities among eigenvalues for $[K|-I]$, those for $[e^{i\gamma}K|-I]$, and those for $[-K|-I]$ in the case that $k_{12}=0$ under the assumption that $f_0=f_N=1$. They are direct consequences  of (iii)--(iv) in Theorem 3.8.\medskip

With the help of  Theorems 3.6--3.7, (3.33), and  a similar method to that used in the proof of Theorem 3.8, one can deduce the following Theorems 3.9--3.11:\medskip

\noindent{\bf Theorem 3.9.} {\it Fix a difference equation $\pmb\omega=(1/f,q,w)$ satisfying that $\prod_{i=0}^{N-1}(1/f_i)>0$. Let  $\gamma\in(-\pi,0)\cup(0,\pi)$ and $K\in SL(2,\mathbb{R})$ satisfy that
$k_{11}-f_0k_{12}\neq0$.
Then the eigenvalues of SLPs $(\pmb\omega, [K|-I])$,
$(\pmb\omega, [e^{i\gamma}K|-I])$, $(\pmb\omega, [-K|-I])$, and $(\pmb\omega, \mathbf{U}_K)$ satisfy the following inequalities:
 \begin{itemize}\vspace{-0.2cm}
\item[{\rm (i)}] if $k_{11}-f_0k_{12}>0$, then
\vspace{-0.2cm} $$\begin{array} {cccc}\lambda_0(K)<\lambda_0(e^{i\gamma}K)<\lambda_0(-K)\leq\lambda_0(\mathbf{U}_K)\leq
\lambda_1(-K)<\\[1.0ex]\lambda_1(e^{i\gamma}K)<\lambda_1(K)
 \leq\lambda_1(\mathbf{U}_K)\leq\cdots\leq\lambda_{N-2}(K)<\lambda_{N-2}(e^{i\gamma}K)<\\[1.0ex]\lambda_{N-2}(-K)
 \leq\lambda_{N-2}(\mathbf{U}_K)\leq\lambda_{N-1}(-K)
 <\lambda_{N-1}(e^{i\gamma}K)<\lambda_{N-1}(K)\end{array}\vspace{-0.2cm}$$
 in the case that $N$ is even;
\vspace{-0.2cm} $$\begin{array} {cccc}\lambda_0(K)<\lambda_0(e^{i\gamma}K)<\lambda_0(-K)\leq\lambda_0(\mathbf{U}_K)\leq
\lambda_1(-K)<\\[1.0ex]\lambda_1(e^{i\gamma}K)<\lambda_1(K)
 \leq\lambda_1(\mathbf{U}_K)\leq\cdots\leq\lambda_{N-2}(-K)<\lambda_{N-2}(e^{i\gamma}K)<\\[1.0ex]\lambda_{N-2}(K)\leq
\lambda_{N-2}(\mathbf{U}_K)\leq\lambda_{N-1}(K)
 <\lambda_{N-1}(e^{i\gamma}K)<\lambda_{N-1}(-K)\end{array}\vspace{-0.2cm}$$
  in the case that $N$ is odd;\vspace{-0.2cm}
  \item[{\rm (ii)}] assertions in  {\rm (i)} hold with $K$ replaced by $-K$.\vspace{-0.2cm}
 \end{itemize}}

\noindent{\bf Theorem 3.10.} {\it Fix a difference equation $\pmb\omega=(1/f,q,w)$ satisfying that $\prod_{i=0}^{N-1}(1/f_i)$ $>0$. Let  $\gamma\in(-\pi,0)\cup(0,\pi)$ and $K\in SL(2,\mathbb{R})$ satisfy that
$k_{11}-f_0k_{12}\neq0$. Then {\rm (i)--(ii)} in Theorem {\rm 3.8}  hold with $k_{12}>0$, $k_{12}<0$, and  $\lambda_n(\mathbf{S}_K)$ replaced by $f_0k_{11}>0$, $f_0k_{11}<0$, and $\lambda_n(\mathbf{T}_K)$,  respectively, where $0\leq n\leq N-1$; {\rm (iii)} in Theorem {\rm 3.8}  holds with $k_{11}>0$, $k_{12}=0$,  and   $\lambda_n(\mathbf{S}_K)$ replaced by $f_0k_{12}<0$, $k_{11}=0$,
 and $\lambda_n(\mathbf{T}_K)$,  respectively, where $0\leq n\leq N-2$; {\rm (iv)} in Theorem {\rm 3.8} also holds. }
\medskip

\noindent{\bf Theorem 3.11.} {\it Fix a difference equation $\pmb\omega=(1/f,q,w)$ satisfying that $\prod_{i=0}^{N-1}(1/f_i)$ $>0$. Let  $\gamma\in(-\pi,0)\cup(0,\pi)$ and $K\in SL(2,\mathbb{R})$ satisfy that
$k_{11}-f_0k_{12}\neq0$. Then {\rm (i)--(ii)} in Theorem {\rm 3.8}  hold with $k_{12}>0$, $k_{12}<0$, and  $\lambda_n(\mathbf{S}_K)$ replaced by $f_0k_{22}-k_{21}>0$, $f_0k_{22}-k_{21}<0$, and $\lambda_n(\mathbf{V}_K)$, respectively, where $0\leq n\leq N-1$; {\rm (iii)} in Theorem {\rm 3.8}  holds with $k_{11}>0$, $k_{12}=0$,  and  $\lambda_n(\mathbf{S}_K)$ replaced by
$k_{11}-f_0k_{12}>0$, $f_0k_{22}-k_{21}=0$, and $\lambda_n(\mathbf{V}_K)$, respectively, where $0\leq n\leq N-2$; {\rm (iv)} in Theorem {\rm 3.8} also holds.  }
\medskip

\noindent{\bf Remark 3.4.} We have not given the similar inequalities as those in  Theorems 3.8--3.11 in the case that $k_{11}-f_0k_{12}=0$ since it is not clear that
what the limits of the polynomial $\Gamma(\lambda)$ given in (3.33) is as $\lambda\rightarrow\pm\infty$ in this case.\medskip

\noindent{\bf 3.4.  Inequalities among eigenvalues for different coupled BCs  }\medskip

In this subsection, we shall  establish inequalities  among eigenvalues for different coupled BCs  applying  Theorems 3.2--3.3.

For each $K\in SL(2,\mathbb{R})$, we set
\vspace{-0.2cm}$$\widehat{K}:=\left(\begin{array}{cc} k_{11} & k_{11}/f_0 \\ k_{21} &(f_0+k_{11}k_{21})/(k_{11}f_0)\end{array}\right)\;\; {\rm if}\;\; k_{11}\neq0;\eqno(3.39)\vspace{-0.2cm}$$
and
\vspace{-0.2cm}$$\widetilde{K}:=\left(\begin{array}{cc} f_0k_{12} & k_{12} \\ (f_0k_{12}k_{22}-1)/k_{12} &k_{22}\end{array}\right) \;\; {\rm if}\;\; k_{12}\neq0.\eqno(3.40)\vspace{-0.2cm}$$
Note that $\widehat{K}, \widetilde{K}\in SL(2,\mathbb{R})$, and  $K=\widehat{K}=\widetilde{K}$ if $k_{11}-f_0k_{12}=0$. The next result
compares  eigenvalues for  $[e^{i\gamma}K|-I]$ with those for $[e^{i\gamma}\widehat{K}|-I]$, and eigenvalues  for $[e^{i\gamma}K|-I]$ with those for $[e^{i\gamma}\widetilde{K}|-I]$, separately. \medskip

\noindent{\bf Theorem 3.12.} {\it Fix a difference equation $\pmb\omega=(1/f,q,w)$. Let $[e^{i\gamma}K|-I]\in\mathcal{B}^\mathbb{C}$, where  $\gamma\in(-\pi,\pi]$ and  $K\in SL(2,\mathbb{R})$ satisfies that $k_{11}-f_0k_{12}\neq0$. Then there are exactly $N$ eigenvalues for $[e^{i\gamma}K|-I]$, and  exactly $N-1$ eigenvalues for both $[e^{i\gamma}\widehat{K}|-I]$ and $[e^{i\gamma}\widetilde{K}|-I]$, where $\widehat{K}$ and $\widetilde{K}$ are defined by {\rm(3.39)--(3.40)}. Furthermore, we have that
\begin{itemize}\vspace{-0.2cm}
\item[{\rm (i)}] if $k_{11}\neq0$, then
\vspace{-0.2cm}$$\begin{array}{cccc}\lambda_0(e^{i\gamma}K)\leq \lambda_0(e^{i\gamma}\widehat{K})\leq \lambda_1(e^{i\gamma}K)\leq \lambda_1(e^{i\gamma}\widehat{K})\leq\cdots\leq\\[1.0ex]\lambda_{N-2}(e^{i\gamma}K)\leq \lambda_{N-2}(e^{i\gamma}\widehat{K})\leq\lambda_{N-1}(e^{i\gamma}K);\end{array}\vspace{-0.2cm}$$
\item[{\rm (ii)}] if $k_{12}\neq0$, then
\vspace{-0.2cm}$$\begin{array}{cccc}\lambda_0(e^{i\gamma}K)\leq \lambda_0(e^{i\gamma}\widetilde{K})\leq \lambda_1(e^{i\gamma}K)\leq \lambda_1(e^{i\gamma}\widetilde{K})\leq\cdots\leq\\[1.0ex]\lambda_{N-2}(e^{i\gamma}K)\leq \lambda_{N-2}(e^{i\gamma}\widetilde{K})\leq\lambda_{N-1}(e^{i\gamma}K).\end{array}\vspace{-0.2cm}$$
\end{itemize}}

\noindent{\bf Proof.}
By Lemma 2.4, the number of    eigenvalues for each BC can be obtained directly.
 Firstly, we show that (i) holds. Let $a_{12}:=k_{12}/k_{11},$ $z:=-e^{i\gamma}/k_{11}$, and $b_{21}:=-k_{21}/k_{11}$.
$\mathbf{A}(s)$ has the same meaning as that in Lemma 2.7.
Then a direct computation implies that $[e^{i\gamma}K|-I]=\mathbf{A}(a_{12})$ and $[e^{i\gamma}\widehat{K}|-I]=\mathbf{A}(1/f_0)$. Hence,  (i) holds by Theorem 3.2.

Assertion (ii) can be shown similarly to that for (i) by Theorem 3.3.  The proof is complete.\medskip

\noindent{\bf Remark 3.5.}  The inequalities in Theorem 3.12 may not be strict. See Example 3.1.\bigskip

\noindent{\bf 4. Inequalities among eigenvalues for different equations }\medskip

In this section,  inequalities among eigenvalues for   equations with different coefficients and weight functions are established
by applying the monotonicity result  of $\lambda_n$ in Theorems 3.1--3.3 in [22].  \medskip

Fix a self-adjoint BC
\vspace{-0.2cm}$$\mathbf{A}=\left [\begin{array} {cccc}a_{11}&a_{12}&b_{11}&b_{12}\\
a_{21}&a_{22}&b_{21}&b_{22}\end{array}  \right ]\vspace{-0.1cm}$$
 in this section.  Let $\mu_1:= a_{11} b_{22}- a_{21} b_{12}$, $\mu_2:= a_{22} b_{12}- a_{12} b_{22}$,  and  $\eta:=-\mu_2/\mu_1$  if $\mu_1\neq0$. If $\mu_1=0$ and $\mu_2=0$, then the BC $\mathbf{A}$ can be written as
$$\vspace{-0.2cm}{\rm either}\begin{array} {llll}  &
 {\mathbf{A}}_1:=\left [\begin{array} {llll} a_{11}&-1&0&0\\
0&0&-1&0\end{array}  \right ] & {\rm or} &
 {\mathbf{A}}_2:=\left [\begin{array} {llll}1& a_{12}&0&0\\
0&0&-1&0\end{array}  \right ].
\end{array}                          \eqno(4.1)         \vspace{-0.2cm}$$

Firstly, we give two lemmas, which play  important roles in establishing inequalities among eigenvalues for equations with different weight functions.
 Fix  $f=\{1/f_n\}_{n=0}^{N}$ and
$q=\{q_n\}_{n=1}^N$. By Lemma 2.4, the number of eigenvalues of  $((1/f,q,w),\mathbf{A})$ is independent of  $w$. Thus, we assume that $((1/f,q,w),\mathbf{A})$ has exactly  $k$ ($N-2\leq k\leq N$) eigenvalues for each $w\in\mathbb{R}^{N,+}$ in the following two lemmas:\medskip

\noindent{\bf Lemma  4.1.} {\it  Fix  $f=\{f_n\}_{n=0}^{N}$, $q=\{q_n\}_{n=1}^N$, $1\leq i\leq N$, $w_1^{(0)},\cdots,w_{i-1}^{(0)}$, $w_{i+1}^{(0)},\cdots,w_{N}^{(0)}$, a self-adjoint BC $\mathbf{A}$, and $1\leq j\leq k$. Let $\lambda_j(w_i):=\lambda_j(w_1^{(0)},\cdots,w_{i-1}^{(0)},w_i,w_{i+1}^{(0)},\cdots,w_{N}^{(0)},\mathbf{A})$ be the $j$-th eigenvalue  function in the $w_i$-direction. If $\lambda_j(w_i^{(0)})=0$ for some $w_i^{(0)}\in\mathbb{R}^+$, then $\lambda_j(w_i)=0$ for all $w_i>w_i^{(0)}$.}\medskip

\noindent{\bf Proof.} Suppose that there exists $w_i'>w_i^{(0)}$ such that $\lambda_j(w_i')\neq0$.  By Theorems 3.1--3.3 in [22],  $\lambda_j(w_i)$ is continuous in $w_i\in\mathbb{R}^{+}$, and
its positive and negative parts  are
non-increasing and  non-decreasing in $w_i$-direction, respectively. This is a contradiction to that $\lambda_j(w_i')\neq0$. The proof is complete.\medskip

\noindent{\bf Lemma  4.2.} {\it Fix  $f=\{f_n\}_{n=0}^{N}$, $q=\{q_n\}_{n=1}^N$, a self-adjoint BC $\mathbf{A}$, and $1\leq j\leq k$. Let $\lambda_j(w):=\lambda_j(1/f,q,w,\mathbf{A})$ be the $j$-th eigenvalue  function for $w\in \mathbb{R}^{N,+}$.  Then
either $\lambda_j(w)\geq0$ for all $w\in\mathbb{R}^{N,+}$ or $\lambda_j(w)\leq0$ for all $w\in\mathbb{R}^{N,+}$.}\medskip

\noindent{\bf Proof. }
In the case that there exists $w^{(0)}\in\mathbb{R}^{N,+}$ such that $\lambda_j(w^{(0)})>0$, we shall show that   $\lambda_j(w_i)\geq0$ for each $1\leq i\leq N$ and  all $w_i\in \mathbb{R}^+$, where $\lambda_j(w_i)$ is defined in Lemma 4.1.    Otherwise, there exists a
$w_i'\in \mathbb{R}^+$ such that $\lambda_j(w_i')<0$. Without loss of generality, assume that  $w_i'<w_i^{(0)}$.
 Then there must exists $w_i''\in(w_i',w_i^{( 0 )})$ such that   $\lambda_j(w_i'')=0$ by the continuity of $\lambda_j(w_i)$ in $w_i\in\mathbb{R}^+$.   By Lemma 4.1, $\lambda_j(w_i)=0$ for all $w_i>w_i''$. This is contradict to $\lambda_j(w_i^{(0)})>0$. Thus  $\lambda_j(w_i)\geq0$ for   all $w_i\in \mathbb{R}^+$. This, together with the monotonicity of $\lambda_j$ in each $w_i$-direction, $1\leq i\leq N$, implies that $\lambda_j(w)\geq0$ for all $w\in\mathbb{R}^{N,+}$.

In the case that there exists $w^{(0)}\in\mathbb{R}^{N,+}$ such that $\lambda_j(w^{(0)})<0$,  with a similar argument above, one can show that  $\lambda_j(w)\leq0$ for all $w\in\mathbb{R}^{N,+}$.

 If it is not one of the above two cases, then $\lambda_j(w)\equiv0$ for all $w\in\mathbb{R}^{N,+}$. The proof is complete.\medskip

Now, inequalities among eigenvalues for equations with different coefficients and weight functions are established. \medskip

\noindent{\bf Theorem  4.1.} {\it Fix a self-adjoint BC ${\mathbf{A}}$.
Consider the following two different  equations: \vspace{-0.2cm}
$$ -\nabla(f_{n}^{(i)}\Delta y_{n})+q_{n}^{(i)}y_{n}=\lambda w_{n}^{(i)}y_{n}, \;\;\;\; \ \  n\in[1,N],\;\;i=1,2,                                                                 \eqno(4.2)_i
\vspace{-0.2cm} $$
and the same BC $\mathbf{A}$.
By $\lambda_n^{(i)} $ denote the $n$-th eigenvalue of ${\rm(4.2)}_{i}$ and $\mathbf{A}$. Let $f_j^{(1)}\leq f_j^{(2)}$ for $0\leq j\leq N-1$, $q_m^{(1)}\leq q_m^{(2)}$ for $1\leq m \leq N$, and $f_N^{(1)}$  and $f_N^{(2)}$ be two given non-zero real  numbers.
  \begin{itemize}\vspace{-0.2cm}
\item[{\rm (i)}]
If one of the following conditions {\rm (1)--(2)} holds,
\begin{itemize}\vspace{-0.2cm}
    \item[{\rm (1)}] $\mu_1\neq0$, $\mu_2\neq0$, and either $ f_0^{(2)}\in(-\infty,1/\eta)$ or $f_0^{(1)}\in(1/\eta,+\infty)$;
 \item[{\rm (2)}] either $\mu_1=0$, $\mu_2\neq0$ or $\mu_1\neq0$, $\mu_2=0$;
\vspace{-0.2cm}\end{itemize}
 then there are exactly $N$ eigenvalues $\lambda_n^{(i)}$ of $(4.2)_i$ and $\mathbf{A}$, where $i=1,2$.  Further, for any given $0\leq n\leq N-1$,
 \begin{itemize}\vspace{-0.2cm}
    \item[{\rm (a)}] if $\lambda_n^{(1)}>0$ and $w_m^{(1)}\geq w_m^{(2)}$, $1\leq m\leq N$,   then
\vspace{-0.2cm}$$\lambda_n^{(1)}\leq\lambda_n^{(2)};\eqno(4.3)\vspace{-0.2cm}$$
   \item[{\rm (b)}] if $\lambda_n^{(1)}\leq0$ and $w_m^{(1)}\leq w_m^{(2)}$, $1\leq m\leq N$,   then
$(4.3)$ holds.
   \vspace{-0.2cm}\end{itemize}
\item[{\rm (ii)}] If one of the following conditions {\rm (3)--(7)} holds,
\begin{itemize}\vspace{-0.2cm}
    \item[{\rm (3)}] $\mu_1\neq0$, $\mu_2\neq0$, and $f_0^{(1)}= f_0^{(2)}=1/\eta$;
 \item[{\rm (4)}] $\mu_1=0$, $\mu_2=0$,  $\mathbf{A}=\mathbf{A}_1$ with $a_{11}\neq0$, and either $ f_0^{(2)}\in(-\infty,-a_{11})$ or $f_0^{(1)}\in(-a_{11},+\infty)$;
 \item[{\rm (5)}] $\mu_1=0$, $\mu_2=0$, $\mathbf{A}=\mathbf{A}_2$ with $a_{12}\neq0$, and either $ f_0^{(2)}\in(-\infty,1/a_{12})$ or $f_0^{(1)}\in(1/a_{12},+\infty)$, where $\mathbf{A}_1$ and $\mathbf{A}_2$ are specified in {\rm(4.1)};
 \item[{\rm (6)}] $\mu_1=0$, $\mu_2=0$,  and $\mathbf{A}=\mathbf{A}_1$ with $a_{11}=0$;
 \item[{\rm (7)}] $\mu_1=0$, $\mu_2=0$,  and $\mathbf{A}=\mathbf{A}_2$ with $a_{12}=0$;
\vspace{-0.2cm}\end{itemize}
 then there are exactly $N-1$ eigenvalues $\lambda_n^{(i)}$ of $(4.2)_i$ and $\mathbf{A}$, where $i=1,2$. Further, for any given $0\leq n \leq N-2$, assertions {\rm(a)--(b)} in {\rm(i)}  hold.
\item[{\rm (iii)}]
If one of the following conditions {\rm (8)--(9)} holds,
\begin{itemize}\vspace{-0.2cm}
    \item[{\rm (8)}] $\mu_1=0$, $\mu_2=0$,  $\mathbf{A}=\mathbf{A}_1$ with $a_{11}\neq0$, and $f_0^{(1)}= f_0^{(2)}=-a_{11}$;
 \item[{\rm (9)}] $\mu_1=0$, $\mu_2=0$, $\mathbf{A}=\mathbf{A}_2$ with $a_{12}\neq0$, and $f_0^{(1)}= f_0^{(2)}=1/a_{12}$;
\vspace{-0.2cm}\end{itemize}
 then there are exactly $N-2$ eigenvalues $\lambda_n^{(i)}$ of $(4.2)_i$ and $\mathbf{A}$, where $i=1,2$. Further, for any given $0\leq n \leq N-3$, assertions {\rm(a)--(b)} in {\rm(i)}  hold. \end{itemize}}

\noindent{\bf Proof.}  The number of eigenvalues of $(4.2)_i$ and $\mathbf{A}$ in each case can be obtained by Lemma 2.4.
Firstly, we show that (i) holds with the assumption (1). Let $0\leq n\leq N-1$.  In the case that $\lambda_n^{(1)}>0$,  $\lambda_n(f^{(1)},q^{(1)},w)\geq 0$ for all $w\in\mathbb{R}^{N,+}$ by Lemma 4.2. By Theorem 3.1 in [22], $\lambda_n(f^{(1)},q^{(1)},w)$ is non-increasing in each $w_m$-direction, $1\leq m\leq N$.  Thus,
 if $w_m^{(1)}\geq w_m^{(2)}$, $1\leq m\leq N$, then
 \vspace{-0.2cm}$$\lambda_n^{(1)}=\lambda_n(f^{(1)},q^{(1)},w^{(1)})\leq \lambda_n(f^{(1)},q^{(1)},w^{(2)}). \eqno(4.4)\vspace{-0.2cm}$$
  Again by Theorem 3.1 in [22], $\lambda_n(f,q,w^{(2)})$ is non-decreasing in $f_j\in(-\infty,1/\eta)$ or $(1/\eta,+\infty)$ in each $f_j$-direction, $0\leq j\leq N-1$; and in $q_m\in\mathbb{R}$ in each $q_m$-direction, $1\leq m\leq N$. Since
  $f_j^{(1)}\leq f_j^{(2)}<1/\eta$ or  $1/\eta<f_j^{(1)}\leq f_j^{(2)}$, $0\leq j\leq N-1$, $q_m^{(1)}\leq q_m^{(2)}$, $1\leq m \leq N$, thus
  \vspace{-0.2cm}$$\lambda_n(f^{(1)},q^{(1)},w^{(2)})\leq\lambda_n(f^{(2)},q^{(2)},w^{(2)})=\lambda_n^{(2)}.  \eqno(4.5)\vspace{-0.2cm}$$
 (4.4)--(4.5) imply (4.3) holds.

   In the case that $\lambda_n^{(1)}\leq0$, with a similar method above, one can show that  (4.3) holds.

  With a similar argument to that in the proof of (i) with the assumption (1), one can show that (i) with the assumption (2), (ii)--(iii) hold. This completes the proof.\medskip

\noindent{\bf Remark  4.1.} Theorem 5.5 of [16] and Theorem 3.6 of [17] give several similar inequalities as those in Theorem 4.1 with the assumption that $f_0^{(1)}=f_0^{(2)}$ and $f_N^{(1)}=f_N^{(2)}$.  In addition, it is required in Theorem 5.5 of [16]   that $w^{(1)}=w^{(2)}$.
Note that it is not required in Theorem 4.1  that $f_N^{(1)}=f_N^{(2)}$ and $w^{(1)}=w^{(2)}$; and it is not required in  (1)--(2) and (4)--(7) in Theorem 4.1 that $f_0^{(1)}=f_0^{(2)}$. Thus, Theorem 4.1  can be regarded as  a generalization of the corresponding results in Theorem 5.5 of [16] and Theorem 3.6 of [17]. \medskip

Combining Theorems 3.2 and  4.1 yields inequalities among eigenvalues of  SLPs with different equations and BCs in $\mathcal{O}_{1,4}^{\mathbb{C}}$. \medskip

\noindent{\bf Corollary  4.1.} Consider the following two different  SLPs:
$(4.2)_i$
and  BCs
\vspace{-0.2cm}$$\mathbf{A}(a_{12}^{(i)},b_{21}^{(i)})=\left [\begin{array} {cccc}1&a_{12}^{(i)}&\bar{z}&0\\
0&z&b_{21}^{(i)}&1\end{array}  \right ],\;\;i=1,2. \eqno(4.6)_i\vspace{-0.2cm}$$
By $\lambda_n^{(i)} $ denote the $n$-th eigenvalue of ${\rm(4.2)}_{i}$ and ${\rm(4.6)}_{i}$. Let $f_j^{(1)}\leq f_j^{(2)}$, $0\leq j\leq N-1$, $q_m^{(1)}\leq q_m^{(2)}$, $1\leq m \leq N$, $f_N^{(1)}$ and $f_N^{(2)}$ be two given non-zero real numbers,  $a_{12}^{(1)}\leq a_{12}^{(2)}$, and $b_{22}^{(1)}\leq b_{22}^{(2)}$.
  \begin{itemize}\vspace{-0.2cm}
\item[{\rm (i)}]
If one of the following two conditions {\rm (1)--(2)} holds,
\begin{itemize}\vspace{-0.2cm}
    \item[{\rm (1)}] $a_{12}^{(1)}\neq0$, $f_0^{(2)}a_{12}^{(1)}>0$,  and either $ a_{12}^{(2)}<1/f_0^{(2)}$ or $f_0^{(1)}>1/a_{12}^{(1)}$;
 \item[{\rm (2)}] $a_{12}^{(1)}=0$, and either $ a_{12}^{(2)}\leq 1/f_0^{(2)}$ or $ 1/f_0^{(2)}<0$;
\vspace{-0.2cm}\end{itemize}
 then there are exactly $N$ eigenvalues $\lambda_n^{(i)}$ of $(4.2)_i$ and $(4.6)_i$, $i=1,2$. Further, for any given $0\leq n\leq N-1$,
 \begin{itemize}\vspace{-0.2cm}
    \item[{\rm (a)}] if $\lambda_n^{(1)}>0$ and $w_m^{(1)}\geq w_m^{(2)}$, $1\leq m\leq N$,   then
\vspace{-0.2cm}$$\lambda_n^{(1)}\leq\lambda_n^{(2)};\eqno(4.7)\vspace{-0.5cm}$$
   \item[{\rm (b)}] if $\lambda_n^{(1)}\leq0$ and $w_m^{(1)}\leq w_m^{(2)}$, $1\leq m\leq N$,   then
$(4.7)$ holds.
   \vspace{-0.2cm}\end{itemize}
\item[{\rm (ii)}] If $a_{12}^{(1)}=a_{12}^{(2)}=1/f_0^{(1)}=1/f_0^{(2)}$, then  there are exactly $N-1$ eigenvalues of $(4.2)_i$ and $(4.6)_i$, $i=1,2$. Further, for any given $0\leq n\leq N-2$, assertions {\rm(a)--(b)} in (i) hold.
   \end{itemize}}

\noindent{\bf Proof.} Firstly, we show that (i) holds with the assumption (1). Direct computations imply that $\mu_1^{(i)}:= a_{11}^{(i)} b_{22}^{(i)}- a_{21}^{(i)} b_{12}^{(i)}=1$, $\mu_2^{(i)}:= a_{22}^{(i)} b_{12}^{(i)}- a_{12}^{(i)} b_{22}^{(i)}=-a_{12}^{(i)}\neq0$,  and  $\eta^{(i)}:=-\mu_2^{(i)}/\mu_1^{(i)}=a_{12}^{(i)}$, $i=1,2$.  If  $a_{12}^{(2)}<1/f_0^{(2)}$,  then $f_0^{(2)}<1/a_{12}^{(1)}=1/\eta^{(1)}$ since $a_{12}^{(1)}\leq a_{12}^{(2)}$ and  $f_0^{(2)}a_{12}^{(1)}>0$.
If $f_0^{(1)}>1/a_{12}^{(1)}$, then $f_0^{(1)}>1/\eta^{(1)}$.
 Fix the BC $\mathbf{A}(a_{12}^{(1)},b_{21}^{(1)})$. By (1) of Theorem 4.1, one gets that there are exactly $N$ eigenvalues of $(4.2)_1$--$(4.6)_1$ and $(4.2)_2$--$(4.6)_1$, and in either case (a) or (b),  for each $0\leq n\leq N-1$,
\vspace{-0.2cm}$$\lambda_n^{(1)}=\lambda_n(1/f^{(1)},q^{(1)},w^{(1)},\mathbf{A}(a_{12}^{(1)},b_{21}^{(1)}))\leq \lambda_n(1/f^{(2)},q^{(2)},w^{(2)},\mathbf{A}(a_{12}^{(1)},b_{21}^{(1)})).\eqno(4.8)\vspace{-0.2cm}$$
Fix the equation $(1/f^{(2)}, q^{(2)}, w^{(2)})$.
If $a_{12}^{(2)}<1/f_0^{(2)}$, then $a_{12}^{(1)}\leq a_{12}^{(2)}<1/f_0^{(2)}$. If $f_0^{(1)}>1/a_{12}^{(1)}$, then $a_{12}^{(2)}\geq a_{12}^{(1)}>1/f_0^{(2)}$.
 By Theorem 3.2, one gets that   there are exactly $N$ eigenvalues of $(4.2)_2$--$(4.6)_1$ and $(4.2)_2$--$(4.6)_2$, and for each $0\leq n\leq N-1$,
\vspace{-0.2cm}$$\lambda_n(1/f^{(2)},q^{(2)},w^{(2)},\mathbf{A}(a_{12}^{(1)},b_{21}^{(1)}))\leq \lambda_n(1/f^{(2)},q^{(2)},w^{(2)},\mathbf{A}(a_{12}^{(2)},b_{21}^{(2)}))=\lambda_n^{(2)}.\eqno(4.9)\vspace{-0.2cm}$$
Combining (4.8)--(4.9) yields that (4.7) holds.

  With a similar argument to that in the proof of (i) with the assumption (1), one can show that  (i) with the assumption (2) and (ii) hold. This completes the proof.\medskip

\noindent{\bf Remark  4.2.}
 One can establish inequalities among eigenvalues of SLPs with different equations and BCs in $\mathcal{O}_{2,4}^{\mathbb{C}}$, $\mathcal{O}_{1,3}^{\mathbb{C}}$, and  $\mathcal{O}_{2,3}^{\mathbb{C}}$, separately, with a similar method to that used in Corollary 4.1. We omit their details.

\bigskip \noindent{\bf \large References}
\def\hang{\hangindent\parindent}
\def\textindent#1{\indent\llap{#1\enspace}\ignorespaces}
\def\re{\par\hang\textindent}
\noindent \vskip 3mm

\re{[1]} F. V. Atkinson,
Discrete and Continuous Boundary Problems, Academic
Press, New York, 1964.

\re{[2]} M. Bohner, O. Dosly, The discrete Pr\"{u}fer transformation, Proc. Amer. Math. Soc. 129 (2001) 2715--2725.

\re{[3]} X. Cao, Q. Kong, H. Wu, A. Zettl, Sturm-Liouville problems whose leading coefficient
function changes sign, Canadian J. Math. 55 (2003) 724--749.

\re{[4]} E. A. Coddington, N. Levinson, Theory of Ordinary Differential Equations, McGraw--Hill, New York, 1955.

\re{[5]} R. Courant, D. Hilbert, Methods of Mathematical Physics, Interscience Publishers, New York, 1953.

\re{[6]} M. Eastham, Q. Kong, H. Wu, A. Zettl,  Inequalities among eigenvalues of Sturm-Liouville problems, J. Inequal. Appl. 3 (1999) 25--43.

\re {[7]} A. Jirari, Second-order Sturm-Liouville difference equations and orthogonal polynomials, Mem. Amer. Math. Soc. 113 (1995).

\re{[8]} Q. Kong, Q. Lin, H. Wu, A. Zettl, A new proof of the inequalities among Sturm-Liouville eigenvalues, PanAmerican Math. J. 10 (2000) 1--11.

\re{[9]} Q. Kong, H. Wu, A. Zettl, Dependence of the $n$-th Sturm-Liouville eigenvalue on the problem, J. Differ. Equ. 156 (1999) 328--354.

\re{[10]} Q. Kong, H. Wu, A. Zettl, Inequalities among eigenvalues of singular Sturm-Liouville problems, Dynamic Systems \& Appl. 8 (1999) 517--531.

\re{[11]} Q. Kong, H. Wu, A. Zettl, Geometric aspects of Sturm-Liouville
problems, I. Structures on spaces of boundary conditions,
Proc. Roy. Soc. Edinb. Sect. A Math. 130 (2000) 561--589.

\re{[12]} Q. Kong, H. Wu, A. Zettl, Left-definite Sturm-Liouville problems, J. Differ. Equ. 177 (2001) 1--26.

\re{[13]} Q. Kong, H. Wu, A. Zettl, Sturm-Liouville problems with finite spectrum, J. Math.
Anal. Appl. 263 (2001) 748--762.

\re{[14]} Q. Kong, H. Wu, A. Zettl, Singular left-definite Sturm-Liouville problems, J. Differ. Equ. 206 (2004) 1--29.

\re{[15]} W. Peng, M. Racovitan, H. Wu, Geometric aspects
of Sturm-Liouville problems, V. Natural loops of boundary conditions
for monotonicity of eigenvalues and their applications, Pac. J. Appl. Math. 4 (2006) 253--273.

\re{[16]} Y. Shi, S. Chen,  Spectral theory of second--order vector
difference equations, J. Math. Anal. Appl. 239
(1999) 195--212.

\re{[17]} Y. Shi, S. Chen,  Spectral theory of higher-order discrete
vector Sturm-Liouville problems, Linear Algebra Appl.
323 (2001) 7--36.

\re{[18]} H. Sun, Y. Shi, Eigenvalues of second-order difference
equations with coupled boundary conditions, Linear Algebra
Appl. 414 (2006) 361--372.

\re{[19]} Y. Wang, Y. Shi,  Eigenvalues of second-order difference
equations with periodic and anti-periodic boundary conditions,
J. Math. Anal. Appl. 309 (2005) 56--69.

\re{[20]} J. Weidmann, Spectral Theory of Ordinary Differential Operators, Lecture Notes in
Mathematics, Vol. 1258, Springer-Verlag, Berlin, 1987.

\re{[21]} A. Zettl, Sturm-Liouville Theory, Mathematical Surveys Monographs, vol. 121, Amer. Math. Soc., 2005.

\re{[22]} H. Zhu, Y. Shi, Continuous dependence of the $n$-th eigenvalue on self-adjoint discrete Sturm-Liouville Problem,  preprint, arXiv:1505.07536 (2015) 1--32.

\re{[23]} H. Zhu, S. Sun, Y. Shi, H. Wu, Dependence of eigenvalues of certain closely discrete Sturm-Liouville  problems, Complex Anal. Oper. Theory, DOI 10.1007/s11785--015--0502--7.

\end{document}